\documentclass[graybox]{svmult}

\usepackage{type1cm}        % activate if the above 3 fonts are
\usepackage{makeidx}         % allows index generation
\usepackage{graphicx}        % standard LaTeX graphics tool
\usepackage{multicol}        % used for the two-column index
\usepackage[bottom]{footmisc}% places footnotes at page bottom

\usepackage{newtxtext}       %
\usepackage{newtxmath}       % selects Times Roman as basic font

\makeindex             % used for the subject index
\begin {document}

\title*{ A multivariate CLT for <<typical>> weighted sums with rate of convergence of order O(1/n) }
\titlerunning{A multivariate CLT for weighted sums}
% Use \titlerunning{Short Title} for an abbreviated version of
% your contribution title if the original one is too long
\author{ Sagak A. Ayvazyan, Vladimir V. Ulyanov}
% Use \authorrunning{Short Title} for an abbreviated version of
% your contribution title if the original one is too long
\institute{Sagak A. Ayvazyan \at Lomonosov Moscow State University, 119991 Moscow, Russia 
%\email{bobkov@math.umn.edu}
%\and Maria A. Danshina \at Moscow Center for Fundamental and Applied Mathematics, Lomonosov Moscow State University, 119991 Moscow, Russia \email{danschina.maria@yandex.ru}
\and Vladimir V. Ulyanov \at 
%Moscow Center for Fundamental and Applied Mathematics, 
Lomonosov Moscow State University, 119991 Moscow, Russia \at National Research University Higher School of Economics, 101000 Moscow, Russia 
\email{vulyanov@cs.msu.ru}}
%
% Use the package "url.sty" to avoid
% problems with special characters
% used in your e-mail or web address
%
\maketitle

\abstract{The "typical" asymptotic behavior of the weighted sums of independent random vectors in $k$-dimensional space is considered. It is shown that in this case the rate of convergence in the multivariate central limit theorem is of order $O(1/n)$. This extends the one-dimensional Klartag and Sodin (2011) result. 
%In this article, we are interested in the analysis of the asymptotic properties of the cycles
%a number of any length in a generalized random graph with random weights of vertices for
%different moment conditions on the distribution of weights.
}

\section {Introduction and Main Result} \label {1}
Let $X, X_1, X_2, \dots, X_n$ be independent identically distributed random vectors in $\mathbb{R}^k$ with finite third absolute moment  $\gamma^3=\mathbb{E}\|X\|^3< \infty$, zero mean $\mathbb{E}X={0}$ and unit covariance matrix $ cov(X)=I$.  
Let $Z$ be the standard Gaussian random variable in $\mathbb{R}^k$ with zero mean and unit covariance matrix.
Denote by $ \mathfrak{B}$, the class of all Borel convex sets in $\mathbb{R}^k$. 

Sazonov~\cite{S68} obtained the following error bound of approximation for distribution of the normalized sum of random vectors by the standard multivariate normal law:  
 \begin{equation}\label{One}\sup\limits_{B \in  \mathfrak{B}}\left|\mathbb{P}\left(\frac{1}{\sqrt n}\sum_{i=1}^{n} X_i \in B \right)-\mathbb{P}(Z \in B) \right|\leq C(k)\frac{\gamma^3}{\sqrt n},
 \end{equation} 
where $C(k)$ depends on dimension $k$ only. 

The bound \eqref{One}  is optimal one in general. Moreover, the rate $O(1/\sqrt{n})$ can not be improved under higher order moment assumptions. This is easy to show in one-dimensional case $k=1$ taking $X$ such that 
$$%\begin{equation}
\mathbb{P}(X=1) = \mathbb{P}(X=-1) = 1/2.\nonumber
$$%\end{equation} 

However, the situation is different when we consider a weighted sum
$$
\theta_1X_1 + \dots + \theta_nX_n,
$$
where $\sum_{j=1}^n\theta_j^2=1$. If we are interested in  
 the typical behavior of these sums for most of $\theta$ in the sense of the normalized Lebesgue measure $\lambda_{n-1}$ on the unit sphere $$S^{n-1}=\{(\theta_1,\dots,\theta_n): \sum_{j=1}^n\theta_j^2=1\},$$ then we have to refer to a recent remarkable result due to Klartag and Sodin. In \cite{KS} they  have showed that in one-dimensional case $k=1$ for any $\rho: 1>\rho>0,$ there exists a set $\mathbb {Q} \subseteq S^{n-1}: \lambda_{n-1}(\mathbb {Q}) > 1 - \rho,$ and a constant $C(\rho)$ depending on $\rho$ only such that for any $\theta=(\theta_1,\dots,\theta_n) \in \mathbb {Q}$ one has
\begin{equation}\label{Two}\sup\limits_{a,b \in \mathbb{R}, a<b}\left|\mathbb{P}\left(a \leq \sum_{i=1}^{n}\theta_iX_i \leq b \right)-\int_a^b \frac{1}{\sqrt{2\pi}}\exp\left(-\frac{x^2}{2}\right)dx \right|\leq C(\rho)\frac{\delta^4}{n},
   \end{equation} 
where $\delta^4=\mathbb{E}\|X\|^4$ and $C(\rho) \leq C \log^2(1/\rho) $ with some absolute constant $C$. It is clear that $C(\rho) \to \infty $ as $\rho$ tends to $0$. And the case of equal weights, that is when $\theta_i = 1/\sqrt{n}$ for all $i = 1, \dots , n$ is the worst case in the sense of closeness of distribution function  of weighted sum to the standard normal distribution function. 

Bobkov \cite{B20} refined the rates of approximation for distributions of weighted sums \eqref{Two} up to order $O(n^{-3/2})$ by using the Edgeworth correction of the fourth order provided $%\beta_5=
\mathbb{E}\|X\|^5 < \infty$. In addition, see \cite{BCG18} and \cite{GNU17} for recent approximation results related to weighted sums in one-dimensional case $k=1.$

 In this paper, we consider the rate of convergence of a weighted sum of independent random vectors to a multivariate standard normal vector. The estimate \eqref{Two} for weighted sums is generalized to the multidimensional case in the form of the following theorem

\begin {theorem} \label {t3}% immediately specify the link
Let $ X_1, X_2,\dots,X_n $ be independent random vectors of dimension $ k $ with zero means $ \mathbb {E} X_j = \overline {0} $, unit covariance matrices $ cov (X_j) = I $ and finite fourth absolute moments $ \delta_j ^ 4 = \mathbb {E} \| X_j \| ^ 4 $ for $ j = 1,\dots,n. $ Denote by $ \mathfrak {B} $   the class of all convex Borel sets and by $ \Phi $  the multidimensional normal distribution with zero mean and unit covariance matrix. Let $ \delta ^ 4 = \frac {1} {n} \sum_ { j = 1} ^ n \delta_j ^ 4 $ and $ \lambda_{n-1} $ be the normalized Lebesgue measure on the unit sphere $ S ^ {n-1} = \{(\theta_1, \theta_2, ,\dots,, \theta_n): \sum_ {j = 1} ^ n \theta_j ^ 2 = 1 \}. $ Then, for any $ \rho> 0 $, there is a subset $ \mathbb {Q} \subseteq S ^ {n-1} $ with $ \lambda_{n-1} (\mathbb {Q})> 1- \rho $ and constant $ C (\rho, k) $ such that for any $ \theta = (\theta_1, ,\dots,, \theta_n) \in \mathbb {Q}, $ one has
    $$ \sup \limits_ {B \in \mathfrak {B}} \Big | P \Big (\sum_ {j = 1} ^ {n} \theta_jX_j \in B \Big) - \Phi (B) \Big | \leq C (\rho, k) \frac {\delta ^ 4} {n}. $$
        Moreover, $ C (\rho, k) \leq C (k) \ln ^ 2 (\frac {1} {\rho}) $, where $ C (k) $ is a universal constant depending only on the dimension $ k. $
\end {theorem}
If we replace the class $ \mathfrak {B} $ with a  smaller class $ \mathfrak {B_0} $ of all centered ellipsoids the situation changes noticeably.
  In this case, the distribution of the normalized sum of i.i.d. random vectors with 
  $ \theta_1 = \dots = \theta_n = 1/\sqrt{n}$ is approximated   by a Gaussian distribution  on the class  $ \mathfrak {B_0} $  with an accuracy of the order from $o(1/\sqrt{n})$ up to $O(1/n)$ under the appropriate dimension of space and when the summands satisfy some moment conditions, for example, finiteness of the fourth absolute moment. See, e.g.
  \cite{12}, \cite{20} and~\cite{14}. For non-i.i.d. random vectors case see~\cite{23}.

\section {Notation and auxiliary results} \label {2}
{\it Notation} In the following, the weight coefficients $ \theta_1, \theta_2 , \dots , \theta_n $, according to the statement of the theorem, will belong to the unit sphere
$$ S ^ {n-1} = \{\theta_1, \theta_2 , \dots , \theta_n: \sum_ {j = 1 } ^ n \theta_j ^ 2 = 1 \}.
$$
Define $ \Theta = (\Theta_1, \Theta_2 ,\dots, \Theta_n) $ as a random vector uniformly distributed on $ S ^ {n-1} $, and for a given set $ \mathbb {Q} \subseteq S ^ {n -1} $ denote the normalized Lebesgue measure on the unit sphere as
\begin {equation} \label {dlambda}
     \lambda_{n-1} (\mathbb {Q}) = P (\Theta \in \mathbb {Q}).\nonumber
  \end {equation}
We denote the absolute fourth-order moment of the random vector $ X_j, \; j = 1,\dots,n, $ as
  \begin {equation} \label {d02}
     \delta_j ^ 4 = \mathbb {E} \| X_j \| ^ 4,
  \end {equation}
  fourth-order weighted absolute moment
\begin {equation} \label {d0}
\delta _ {\theta} ^ 4 = \sum_ {j = 1} ^ {n} \theta ^ 4_j \delta ^ 4_j,
  \end {equation}
and the averaged absolute moment of the fourth order
  \begin {equation} \label {d01}
   \delta ^ 4 = \frac {1} {n} \sum_ {j = 1} ^ n \delta_j ^ 4.
    \end {equation}
We introduce truncated random variables
        $ Y_j $ and $ Z_j, j = 1,\dots,n, $ as
        \begin {equation} \label {d2}
            Y_j = X_j 1_X (\| \theta_j X_j \| \leq 1),
        \end {equation}
        \begin {equation} \label {d3}
           Z_j = Y_j- \mathbb {E} Y_j,\nonumber
        \end {equation}
        where $ 1_X(A) $ is the indicator function of an event $A$.
        For these random variables, we define the weighted expectation and the weighted covariance matrix as
        \begin {equation} \label {d4}
           A_n = \sum_ {j = 1} ^ {n} \theta_j \mathbb {E} Y_j,\nonumber
        \end {equation}
        \begin {equation} \label {d5}
           D = \sum_ {j = 1} ^ {n} \theta_j ^ 2cov (Z_j),
        \end {equation}
        \begin {equation} \label {d6}
          Q ^ 2 = D ^ {- 1}.
        \end {equation}
        We also use the notation for the distributions of random vectors appeared  in the proof: over all Borel sets $ B $
        \begin {equation} \label {df1}
            F_ {X} (B) = P \Big (\sum_ {j = 1} ^ n \theta_jX_j \in B \Big),
        \end {equation}
        \begin {equation} \label {df2}
            F_ {Y} (B) = P \Big (\sum_ {j = 1} ^ n \theta_jY_j \in B \Big),
        \end {equation}
        \begin {equation} \label {df3}
            F_ {Z} (B) = P \Big (\sum_ {j = 1} ^ n \theta_jZ_j \in B \Big),\nonumber
        \end {equation}
        \begin {equation} \label {df4}
            F_ {X_j} (B) = P (\theta_jX_j \in B), \; j = 1,\dots,n,
        \end {equation}
                \begin {equation} \label {df4_1}
            F_ {Y_j} (B) = P (\theta_jY_j \in B), \; j = 1,\dots,n.
        \end {equation}
 
Also $ \Phi_ {a, V} $ will denote the distribution of the multidimensional standard normal law with expectation $ a $ and covariance matrix $ V. $
The proof will use the technique of characteristic functions. The characteristic function of the random vector $ Z_j $ is defined as
\begin {equation} \label {df5}
\varphi_j (t) = \mathbb {E} \exp (i tZ_j), \; j = 1,\dots,n.
\end {equation}
We denote $ \hat {F} _Z = \prod_ {j = 1} ^ {n} \varphi_j (\theta_jt) $ and $ \hat {\Phi} _ {a, V} $ as the corresponding characteristic functions of weighted sums of random vectors $ \sum_ {j = 1} ^ {n} Z_j \theta_j, $ and a random vector with multivariate normal distribution.
 
The absolute moment of order $ s $ of the random vector $ Z $ is defined as
$ \rho_s (Z) = \mathbb {E} \| Z \| ^ s. $
Also, for the coefficients $ \theta_1, \theta_2 ,\dots, \theta_n $, the weighted absolute moment of order $ s $ is determined as
      \begin {equation} \label {d7}
           \rho_s = \sum_ {j = 1} ^ n \rho_s (\theta_jZ_j)
        \end {equation}
and
     \begin {equation} \label {d8}
           \eta_s = \sum_ {j = 1} ^ n \rho_s (Q \theta_jZ_j),
        \end {equation}
for $ 1 \leq m \leq s-1 $ we have
     $$
            \rho_s (Q \theta_jZ_j) \leq \| Q \| ^ {s} \rho_s (\theta_jZ_j))
            = \| Q \| ^ {s} \rho_s (\theta_j (Y_j- \mathbb {E} Y_j))
     $$
              \begin {equation} \label {d8b}
            \leq \| Q \| ^ {s} 2 ^ s \rho_s (\theta_jY_j)
            \leq \| Q \| ^ {s} 2 ^ s \rho_ {s-m} (\theta_jX_j).
              \end {equation}
For a given nonnegative vector $ \alpha $, we define $ \alpha $ the moment of the random vector $ Z $ as
\begin {equation} \label {d9.1}
   \mu_ \alpha (Z) = \mathbb {E} Z ^ \alpha,
\end {equation}
for this value the following inequality holds
     \begin {equation} \label {d9}
           \mu_ \alpha (Z) \leq \rho_ {| \alpha |} (Z).
        \end {equation}
Let the random vector $ Z_j $ have finite absolute moments of order $ m. $ Then the characteristic function in a neighborhood of zero satisfies the Taylor expansion
\begin{equation}\label{cumul}
\varphi_j (t) = 1 + \sum_ {1 \leq | \nu | \leq m} \mu_ \nu (Z_j) \frac {(it) ^ \nu} {\nu!} + o (\| t \| ^ m),\;j=1\dots n, \nonumber
\end{equation}
as $ t $ tends to zero.
Next, we define the logarithm of a nonzero complex number as $ z = r \exp (i \xi) $ as
$$ \log (z) = \log (r) + i \xi, $$
where $$ r> 0, \; \; \xi \in (- \pi, \pi]. $$
Thus, we always take the so-called main branch of the logarithm. The characteristic function of the random vector $ Z_j $ is continuous and equal to one at zero. Consequently, in a neighborhood of zero, the Taylor expansion takes place
$$ \log (\varphi_j (t)) = \sum_ {1 \leq | \nu | \leq m} \kappa_ \nu (Z_j) \frac {(it) ^ \nu} {\nu!} + o (\| t \| ^ m), \;j=1\dots n.$$
as $ t $ tends to zero. The expansion coefficients of the logarithm of the characteristic function $ \kappa_ \nu $ are called the cumulants of the random vector $ Z_j. $ The cumulants $ \kappa_\nu $  are explicitly expressed in terms of moments (see Ch.2 Sect.6 in \cite{B}). In particular, the following inequality holds:
     \begin {equation} \label {d10}
           | \kappa_ \nu (Z_j) | \leq c_1 (\nu) \rho_ {| \nu |} (Z_j),\;j=1\dots n.
        \end {equation}
We also point out the most important property of the cumulants of the sum of random vectors and denote the cumulants of the weighted sum as 
     \begin {equation} \label {d11}
          \kappa_ \nu= \kappa_ \nu \Big(\sum_ {j = 1} ^ n \theta_jZ_j\Big) = \sum_ {j = 1} ^ n \theta_j ^ {| \nu |} \kappa_ \nu (Z_j)
        \end {equation}
and
     \begin {equation} \label {d12}
           \kappa_r (t)=\kappa_r \Big(\sum_ {j = 1 } ^ n \theta_jZ_j, t\Big) = \sum_ {| \nu | = r} \frac {\kappa_ \nu  t ^ \nu} {\nu!}.
        \end {equation}
   Also for the characteristic function of the weighted sum, in a neighborhood of zero, one has
$$ \log \Big(\prod_{j=1}^n\varphi_j (t)\Big) = \sum_ {r = 1} ^ m \frac {\kappa_r ( it)} {r!} + o (\| t \| ^ m),$$
let us define the polynomials $ P_r (t,\kappa_ \nu) $ from the formal expression
\begin {equation} \label {d13}
1+ \sum_ {r = 1} ^ \infty P_r (t, \kappa_ \nu) u ^ r = \exp \Big (\sum_ {s = 1} ^ \infty \frac {\kappa_ {s +2} (t)} {(s + 2)!} u ^ s \Big),
\end {equation}
explicitly, we obtain the expressions
$$ P_0 (t, \kappa_ \nu) = 1, $$
\begin {equation} \label {d14}
P_r (t, \kappa_ \nu ) = \sum_{m = 1}^r \frac{1}{m!} \sum_{i_1,\dots,i_m: \sum_{j = 1} ^ mi_j = r } \Bigg[\sum_{\nu_1 ,\dots, \nu_m: | \nu_j | = i_j + 2} \frac {\kappa_{\nu_1} \dots\kappa_{\nu_m}} {\nu_1! \dots \nu_m!} \Bigg] t^{\nu_1 + \dots + \nu_m}.
\end {equation}
For a given positive vector $ \alpha $ and a number $ m $, we denote $ D ^ \alpha $ and $ D_m $ as differential operators
$$ D ^ \alpha f (t) = \frac {\partial ^ {\alpha_1 + \alpha_2 + \dots + \alpha_k} f (t)} {(\partial t_1) ^ {\alpha_1} (\partial t_2 ) ^ {\alpha_2} \dots (\partial t_k) ^ {\alpha_k}} $$
and
  $$ D_m f (t) = \frac {\partial f (t)} {\partial t_m}. $$
  \begin {lemma} \label {lemma1}
  Let $ Z_1, Z_2,\dots,Z_n $ be independent random vectors (nondegenerate at zero) with a finite absolute moment of order $ s. $ Then for any $ 2 <r \leq s $
  $$ \Big (\frac {\rho_r} {\rho_2 ^ {\frac {r} {2}}} \Big) ^ {\frac {1} {r-2}} \leq \Big (\frac { \rho_s} {\rho_2 ^ {\frac {s} {2}}} \Big) ^ {\frac {1} {s-2}}, $$
  where $ \rho_s $ is defined in (\ref {d7}).
  \end {lemma}
  \begin {proof}
  The proof of Lemma follows the scheme of the proof of Lemma 6.2 \cite {B}.
Let us show that $ \log \rho_r $ is a convex function on $ [2, s]. $, where $\rho_s = \sum_ {j = 1} ^ n \rho_s (\theta_jZ_j)= \sum_ {j = 1} ^ n\mathbb {E} \|\theta_j Z_j \| ^ s$ ( see \ref{d7}). This follows from Holder's inequality for $ \alpha + \beta = 1 $ and $ r_1, r_2 \in [2, s] $
$$ \rho _ {\alpha r_1 + \beta r_2} (Z_j) \leq \rho_ {r_1} ^ \alpha (Z_j) \rho_ {r_2} ^ \beta (Z_j), $$
$$ \rho _ {\alpha r_1 + \beta r_2} \leq \sum_ {j = 1} ^ n \rho_ {r_1} ^ \alpha (\theta_jZ_j) \rho_ {r_2} ^ \beta (\theta_jZ_j) $$
$$ \leq \Big (\sum_ {j = 1} ^ n \rho_ {r_1} (\theta_jZ_j) \Big) ^ \alpha \Big (\sum_ {j = 1} ^ n \rho_ {r_2} (\theta_jZ_j) \Big) ^ \beta = \rho_ {r_1} ^ \alpha \rho_ {r_2} ^ \beta. $$
so
$$ \log (\rho _ {\alpha r_1 + \beta r_2}) \leq \alpha \log (\rho_ {r_1}) + \beta \log (\rho_ {r_2}).
$$
Further, suppose that $ \rho_2 = 1 $, then
$$ \log (\rho_r ^ {\frac {1} {r-2}}) = \frac {\log (\rho_r) - \log (\rho_2)} {r-2} $$
increases due to the fact that it is the slope between the points $ (2, \rho_2) $ and $ (r, \rho_r) $ of the function $ \rho_r. $ In the general case, when $ \rho_2 \neq 1 $, it is necessary to consider the random vectors $ \hat {Z} _j = Z_j/ \sqrt \rho_2, \; \; j = 1,\dots,n. $
 
  \end {proof}
    \begin {lemma} \label {lemma2}
Let $ X_1, X_2,\dots,X_n $ be independent random vectors with zero mean, unit covariance matrix and finite fourth absolute moment. If the following condition holds $\delta_ \theta ^ 4 \leq (8k)^{-1}, $
then the weighted covariance matrix $ D $ satisfies
$$ \Big | \langle t, Dt \rangle - \| t \| ^ 2 \Big | \leq 2k \delta_ \theta ^ 4 \| t \| ^ 2 $$
and
$$ \| D-I \| \leq \frac {1} {4}, \; \; \; \frac {3} {4} \leq \| D \| \leq \frac {5} {4},\;\;\;\| D ^ {- 1} \| \leq \frac {4} {3}. $$
Where $ \delta_ \theta ^ 4 $ and matrix $ D $ are defined in (\ref {d0}) and (\ref {d5}).
  \end {lemma}
\begin {proof}
The proof of Lemma follows the scheme of the proof of Corollary 14.2 \cite {B}.
First, we prove two auxiliary inequalities for the mathematical expectations of the original and truncated random vectors
$$ \Big | \mathbb {E} \theta_jY_ {ji} \Big | = \Big | \mathbb {E} X_ {ji} \theta_j1_X (\| X_j \theta_j \|> 1) \Big |$$
$$\leq \mathbb {E} \| X_j \theta_j \| 1_X (\| X_j \theta_j \|> 1) =\theta_j ^4\mathbb {E} \| X_j\|^4 \leq \theta_j ^ 4 \delta_j ^ 4, $$
and
$$ \theta_j ^ 2 \Big | \mathbb {E} X_ {ji} X_ {jl} - \mathbb {E} Y_ {ji} Y_ {jl} \Big | = \theta_j ^ 2 \Big | \mathbb {E} X_ {ji} X_ {jl} 1_X (\| \theta_jX_ {ji} \|> 1) \Big | $$
$$ \leq \mathbb {E} \| \theta_j X_j \| ^ 2 1_X (\| \theta_j X_ {ji} \|> 1) \leq \theta_j ^ 4 \mathbb {E} \|  X_j \| ^ 4 = \theta_j ^ 4 \delta_j ^ 4. $$
Also, note that
$$ | \mathbb {E} \theta_jY_ {ji} | = | \mathbb {E} \theta_jX_ {ji} 1_X (\| X_j \theta_j \| \leq 1) | \leq 1. $$
Define  Kronecker delta function  as $\delta_{ij}=1_X(i=j),$ for $i,j=1,\dots,k.$
 By definition of covariance matrix  of the weighted sums  $D = \sum_ {j = 1} ^ {n} \theta_j ^ 2cov (Z_j)$ ( see (\ref{d5})), one has
$$ \Big | \langle t, Dt \rangle - \langle t, t \rangle \Big | = \Big | \sum_ {i, l} ^ kt_ {i} t_l (d_ {il} -\delta_ {il}) \Big |, $$
wherein
$$ \Big | d_ {il} -\delta_ {il} \Big | \leq \sum_ {j = 1} ^ n \theta_j ^ 2 \Big | cov (X_ {ji}, X_ {jl}) - cov (Y_ {ji}, Y_ {jl}) \Big | $$
$$ \leq \sum_ {j = 1} ^ n \theta_j ^ 2 \Big | \mathbb {E} X_ {ji} X_ {jl} - \mathbb {E} Y_ {ji} Y_ {jl} + \mathbb {E} Y_ {ji} \mathbb {E} Y_ {jl} \Big | \leq \sum_ {j = 1} ^ n \Big (\theta_j ^ 4 \delta_j ^ 4 + \theta_j ^ 4 \delta_j ^ 4 \Big) = 2 \delta_ \theta ^ 4, $$
where  $\delta_j ^ 4 = \mathbb {E} \| X_j \| ^ 4, \delta ^ 4 = \frac {1} {n} \sum_ {j = 1} ^ n \delta_j ^ 4$ ( see (\ref{d0}), (\ref{d01})).
Finally we get
$$ \Big | \sum_ {i, l} ^ kt_it_l (d_ {il} -\delta_{il}) \Big | \leq2 \delta_ \theta ^ 4 \Big (\sum_i ^ k | t_i | \Big) ^ 2 \leq2k \delta_ \theta ^ 4 \| t \| ^ 2. $$
By the definition of the matrix norm, it follows that
  $$ \| D-I \| = \sup \limits _ {\| t \| \leq 1} \Big | \langle t, (D-I) t \rangle \Big | \leq2k \delta_ \theta ^ 4. $$
Since  $ \delta_ \theta ^ 4 <(8k) ^ {- 1} $, then
  $$ \| D-I \| \leq \frac {1} {4},\;\;\;\frac {3} {4} \leq \| D \| \leq \frac {5} {4}. $$
Further,
  $$ \langle t, Dt \rangle \geq \| t \| ^ 2- \frac {1} {4} \| t \| ^ 2 = \frac {3} {4} \| t \| ^ 2. $$
  Therefore, for the 
  inverse matrix $D^{- 1}$    one has
  $$ \| D ^ {- 1} \| \leq \frac {4} {3}. $$
  \end {proof}
   \begin {lemma} \label {lemma3}
   Let $ Z_1, Z_2,\dots,Z_n $ be independent random vectors with zero mean, nondegenerate covariance matrix $ D $ (see (\ref {d5})) and finite absolute moment $ \rho_s $ of order $ s \geq 3. $ Then there are constants $ c_1 (k, s) $
   and $ c_2 (k, s) $ such that for any $ | \alpha | \leq s $ and $ \| t \| \leq c_1 (k, s) \min \Big\{\eta_s ^ {- \frac {1} {s}}, \eta_s ^ {- \frac {1} {s-2}} \Big\} $ one has
  $$ \Big | D ^ \alpha \Big [\prod_ {j = 1} ^ {n} \varphi_j (\theta_j Qt) - \exp \Big (- \frac {1} {2} \| t \| ^ 2 \Big) \sum_ {r = 0} ^ {s-3} P_r (iQt, \kappa_ \nu) \Big] \Big | $$
  $$ \leq c_2 (k, s) \eta_s (\| t \| ^ {s- | \alpha |} + \| t \| ^ {3 (s-2) + | \alpha |}) \exp \Big (- \frac {1} {4} \| t \| ^ 2 \Big), $$
  where $ \varphi (t) $ is defined in (\ref {df5}), matrix $ Q $ in (\ref {d6}), $ \eta_s $ in (\ref {d8}) and polynomial $ P_r (t, \kappa_ \nu) $ to (\ref {d13}), (\ref {d14}).
   \end {lemma}
    \begin {proof}
    The proof follows the scheme of the proof of Theorem 9.11 \cite {B}. First, suppose that the covariance matrix is $ D = I $, then from the definition of weighted absolute moments $\rho_s = \sum_ {j = 1} ^ n \rho_s (\theta_jZ_j),\;\eta_s = \sum_ {j = 1} ^ n \rho_s (Q \theta_jZ_j)$( see (\ref{d7})(\ref{d8})), where $ Q ^ 2 = D ^ {- 1}$ ( see (\ref{d6})), it follows that $ \eta_s = \rho_s.$ Using Holder's inequality, we obtain that for the characteristic function of the random vector $ Z_j $ at $ \| t \| \leq \rho_s ^ {- \frac {1} {s}} $ the following inequality holds $$| \varphi_j (\theta_j t) -1 | \leq \frac {\mathbb {E} \langle t, \theta_
    jZ_j \rangle ^ 2} {2}  \leq \frac {\| t \| ^ 2\mathbb {E} \| \theta_jZ_j \| ^ 2} {2}  \leq \frac {\| t \| ^ 2} {2} \Big (\mathbb {E} \| \theta_jZ_j \| ^ s \Big) ^ {\frac {2} {s}} \leq \frac {\| t \| ^ 2} {2} \rho_s ^ {\frac {2} {s}} \leq \frac {1} {2}, $$
    therefore, the characteristic function of the random vector $ \theta_jZ_j $, defined as $\varphi_j (t) = \mathbb {E} \exp (i tZ_j) $( see \ref{df5}), does not vanish in a given interval, and therefore the following functions can be defined
    $$ h_j (t) = \log (\varphi_j (\theta_jt)) - \Big (- \theta_j ^ 2 \| t \| ^ 2 + \sum_ {r = 1} ^ {s-3} \frac { \kappa_ {r + 2} (Z_j, it)} {(r + 2)!} \Big), $$
    $$ h (t) = \sum_j ^ nh_j (t), $$
    $$ \zeta (t) = \frac {1} {2} \| t \| ^ 2 + \sum_ {r = 1} ^ {s-3} \frac {\kappa_ {r + 2} (it)} {(r + 2)!}, $$
    where $\kappa_r (t)= \sum_ {| \nu | = r} \frac {\kappa_ \nu  t ^ \nu} {\nu!}$ ( see (\ref{d11}), (\ref{d12})).
     Since the weighted absolute moment of the second order is $ \rho_2 = k $, then by Lemma \ref {lemma1} for $ 2 <r \leq s $,
    $$ (\rho_r) ^ {\frac {1} {r-2}} \leq (\rho_s) ^ {\frac {1} {s-2}} (\rho_2) ^ {\frac {1} { r-2} - \frac {1} {s-2}} = (\rho_s) ^ {\frac {1} {s-2}} (k) ^ {\frac {1} {r-2} - \frac {1} {s-2}} \leq (\rho_s) ^ {\frac {1} {s-2}} k $$
    and for the cumulants of the distribution $ \kappa_ \nu, $ for $ 2 <| \nu | \leq s $, due to (\ref {d9}), (\ref {d10})
    $$ | \kappa_ \nu | ^ {\frac {1} {| \nu | -2}} \leq (c_1 (\nu) \rho_ {| \nu |}) ^ {\frac {1} {| \nu | -2}} \leq c (\nu) k \rho_s ^ {\frac {1} {s-2}} \leq \hat {c} _1 (s, k) \rho_s ^ {\frac {1} { s-2}}. $$
Next, consider the expression
     $$ D ^ \alpha \Big [\prod_ {j = 1} ^ n \varphi_j (\theta_jt) - \exp (\zeta (t)) \Big] = D ^ \alpha \Big [(\exp (h (t) -1) \exp (\zeta (t)) \Big] $$
     \begin {equation} \label {e1}
    = \sum \limits_ {0 \leq \beta \leq a} c (\alpha, \beta) D ^ \beta \exp (\zeta (t)) D ^ {\alpha- \beta} (\exp (h (t)) - 1).
    \end {equation}
    Denote $$ c (s, k) = \sum \limits_ {r = 0} ^ {s-3} \sum \limits_ {| \nu | = r + 2} \frac {1} {\nu!} $$
and  notice that for $ \| t \| \leq \Big (\hat {c} _1 (s, k) 8c (s, k) \rho_s ^ {\frac {1} {s-2}} \Big) ^ {- 1} $ the following chain of inequalities holds
    $$ \Big | \sum_ {r = 0} ^ {s-3} \sum_ {| \nu | = r + 2} \frac {\kappa_ \nu} {\nu!} t ^ \nu \Big | \leq \sum_ {r = 0} ^ {s-3} \sum_ {| \nu | = r + 2} \| t \| ^ {r + 2} \frac {(\hat {c} _1 (s, k) \rho_s ^ {\frac {1} {s-2}}) ^ r} {\nu!} $$
$$ \leq \| t \| ^ 2 \sum_ {r = 0} ^ {s-3} \sum_ {| \nu | = r + 2} \frac {(\hat {c} _1 (s, k) \| t \| \rho_s^{\frac{1}{s-2}}) ^ r} {\nu!} $$
$$ \leq \| t \| ^ 2 \sum_ {r = 0} ^ {s-3} \sum_ {| \nu | = r + 2} \frac {1} {\nu!} \Big (\hat {c} _1 (s, k) \rho_s ^ {\frac {1} {s-2}} \rho_s ^ {- \frac {1} {s-2}} (\hat {c} _1 (s, k) 8c (s, k)) ^ {- 1} \Big) ^ r $$
\begin {equation} \label {e11}
     \leq \| t \| ^ 2 \sum_ {r = 0} ^ {s-3} \sum_ {| \nu | = r + 2} \frac {1} {\nu! 8c (s, k) ^ r} \leq \frac {\| t \| ^ 2} {8c (s, k)} \sum_ {r = 0} ^ {s-3} \sum_ {| \nu | = r + 2} \frac {1} {\nu!} = \frac {\| t \| ^ 2} {8},
\end {equation}    
therefore, the module of the function $ \zeta (t) $ is bounded
$$ | \zeta (t) | \leq \frac {\| t \| ^ 2} {2} + \frac {\| t \| ^ 2} {8} = \frac {5 \| t \| ^ 2} {8}. $$
Similarly, it can be shown that the modulus of the derivative of this function
$$
\Big | D ^ \beta \zeta (t) \Big | = \Big | D ^ \beta \sum_ {r = 0} ^ {s-3} \sum_ {| \nu | = r + 2} \frac {\kappa_ \nu} {\nu!} t ^ \nu \Big | 
= \Big | \sum_ {r = \max \{0, | \beta | -2 \}} ^ {s-3} \sum_ {| \nu | = r + 2, \nu \geq \beta} \frac {\kappa_ \nu} {(\nu- \beta)!} t ^ {\nu- \beta} \Big |
$$
$$ \leq \sum_ {r = \max \{0, | \beta | -2 \}} ^ {s-3} \sum_ {| \nu | = r + 2, \nu \geq \beta} \frac {(\hat {c} _1 (s, k) \rho_s) ^ {\frac {r} {s-2}}} {(\nu- \beta)!} \| t \| ^ {r + 2- | \beta |} $$
$$
\leq \sum_ {r = \max \{0, | \beta | -2 \}} ^ {s-3} \sum_ {| \nu | = r + 2, \nu \geq \beta} \frac {(\hat {c} _1 (s, k) \rho_s) ^ {\frac {r} {s-2}}} {(\nu- \beta)!} \Big (\hat {c} _1 (s, k) 8c (s, k) \rho_s ^ {\frac {1} {s-2}} \Big) ^ {- r} \| t \| ^ {2- | \beta |}
$$
$$
\leq c_2 (s, k, \beta) \| t \| ^ {2- | \beta |}.
$$
Further, let $ j_1, j_2,\dots,j_r $
be non-negative numbers 
$ \beta_1, \beta_2 ,\dots, \beta_r$ 
non-negative vectors satisfying the equality
$ \sum_ {i = 1} ^ rj_i \beta_i = \beta, $
since $$ \| t \| ^ {2 \sum \limits_ {i = 1} ^ rj_i} \leq \| t \| ^ {2} + \| t \| ^ {2 -| \beta |}, 
$$
then the derivative $ \zeta (t) $ has the following representation for $ t \neq0 $
$$ \Big | (D ^ {\beta_1} \zeta (t)) ^ {j_1} \dots (D ^ {\beta_r} \zeta (t)) ^ {j_r} \Big | \leq c_ {3} (s, k) \| t \| ^ {\sum \limits_ {i = 1} ^ rj_i (2- | \beta_j |)} $$
$$ \leq c_ {3} (s, k) (\| t \| ^ {2- | \beta |} + \| t \| ^ {| \beta |}), $$
Lemma 9.2 \cite {B} implies that for $ \| t \| \leq \Big (\hat {c} _1 (s, k) 8c (s, k) \rho_s ^ {\frac {1} {s-2}} \Big) ^ {- 1}$
$$ \Big | D ^ \beta \exp (\zeta (t)) \Big | \leq c_ {4} (s, k) (\| t \| ^ {2- | \beta |} + \| t \| ^ {| \beta |}) \exp (\zeta (t)) $$
\begin {equation} \label {e2}
    \leq c_ {4} (s, k) (\| t \| ^ {2- | \beta |} + \| t \| ^ {| \beta |}) \exp \Big (- \frac {3 } {8} \| t \| ^ 2 \Big).
    \end {equation}
Further, for $ \beta: \; 0 \leq | \beta | \leq s, $ one has 
$$ D ^ {\hat {\beta}} (D ^ \beta h_j) (0) = 0, \; [0 \leq | \hat {\beta} | \leq s- | \beta | -1] . $$
Therefore, applying Corollary 8.3 \cite {B} and the fact that $g \equiv D ^ \beta h_j $
we obtain the inequality
$$ \Big | D ^ \beta h_j (t) \Big | \leq \sum \limits_ {| \hat {\beta} | = s- | \beta |} \frac {| t ^ {\hat {\beta}} |} {\hat {\beta}!} \sup \Big (\Big | (D ^ {\hat {\beta}} g) (ut) \Big |: 0 \leq u \leq 1 \Big). $$
If $ | \hat {\beta} | = s- | \beta | $, then, by Lemma 9.4 \cite {B}, the following relation holds:
$$ \Big | D ^ {\hat {\beta}} g (ut) \Big | = \Big | D ^ {\hat {\beta} + \beta} h_j (ut) \Big | = \Big | D ^ {\hat {\beta} + \beta} \log \varphi_j (\theta_jut) \Big || \theta_j | ^ s\leq | \theta_j | ^ sc_2 (s) \rho_ {s} (Z_j), $$
so
\begin {equation} \label {e3}
    \Big | D ^ \beta h (t) \Big | \leq c_2 (s) \rho_s \| t \| ^ {s- | \beta |} \Big (\sum \limits_ {| \hat {\beta} | = s- | \beta |} \frac {1} {\hat {\beta}} \Big).
\end {equation}
If $ \beta = 0 $, then similarly as in (\ref {e11}) we get that
$$ | h (t) | \leq c_2 (s) \Big (\sum_ {| \hat {\beta} | = s} \frac {1} {\hat {\beta}} \Big) \rho_s \| t \| ^ s \leq \frac {\| t \| ^ 2} {8}. $$
If $ \alpha- \beta = 0 $
$$ \Big | D ^ {\alpha- \beta} (\exp (h (t)) - 1)) \Big | = \Big | \exp (h (t)) - 1 \Big | \leq | h (t) | \exp (| h (t) |) $$
\begin {equation} \label {e4}
    \leq c_ {5} (s, k) \rho_s \| t \| ^ s \exp \Big (\frac {\| t \| ^ 2} {8} \Big).
\end {equation}
If $ \alpha > \beta $
$$ D ^ {\alpha- \beta} \Big (\exp (h (t)) - 1 \Big) = D ^ {\alpha- \beta} \exp (h (t)), $$
then in this case the derivative is represented as a linear combination of the following form
$$ (D ^ {\beta_1} h (t)) ^ {j_1} \dots (D ^ {\beta_r} h (t)) ^ {j_r} \exp (h (t)), $$
where $ \sum_ {i = 1} ^ rj_i \beta_i = \alpha- \beta. $ \\
From (\ref {e3}) and the inequality $ \| x \| ^ a \leq \| x \| ^ b + \| x \| ^ c, \; \; 0 \leq b \leq a \leq c $, it follows that
$$ \Big | (D ^ {\beta_1} h (t)) ^ {j_1} \dots (D ^ {\beta_r} h (t)) ^ {j_r} \Big |  \leq c_ {6} (s, k) \rho_s \Big (\| t \| \rho_s ^ {\frac {1} {s-2}} \Big) ^ {(s-2) (\sum \limits_ {i = 1} ^ rj_i -1)} \| t \| ^ {s-2 + 2 \sum \limits_ {i = 1} ^ rj_i- | \alpha- \beta |} $$
$$ \leq c_ {7} (s, k) \rho_s (\| t \| ^ {s- | \alpha- \beta |} + \| t \| ^ {s- | \alpha- \beta | -2}). $$
Therefore, if $ \alpha> \beta $, then
\begin {equation} \label {e5}
    \Big | D ^ {\alpha- \beta} (\exp (h (t)) - 1) \Big | \leq c_ {8} (s, k) \rho_s (\| t \| ^ {s- | \alpha- \beta |} + \| t \| ^ {s- | \alpha- \beta | +2 }) \exp \Big (\frac {\| t \| ^ 2} {8} \Big).
\end {equation}
Using (\ref {e2}), (\ref {e4}), (\ref {e5}) in (\ref {e1}), we get
$$ D ^ {\alpha} \Big [\prod_ {j = 1} ^ n \varphi_j (\theta_jt) - \exp \Big (- \frac {\| t \| ^ 2} {2} + \sum_ {r = 1} ^ {s-3} \frac {\kappa_ {r + 2} (it)} {(r + 2)!} \Big) \Big] $$
$$ \leq c_ {9} (s, k) \rho_s (\| t \| ^ {s- | \alpha|} + \| t \| ^ {s- | \alpha | -2}) \exp \Big (- \frac {\| t \| ^ 2} {4} \Big). $$
Next, we use Lemma 9.7 \cite {B} when setting $ u = 1 $ and the inequality
    $$ | \kappa_ \nu | ^ {\frac {1} {| \nu | -2}} \leq \hat {c} _1 (s, k) \rho_s ^ {\frac {1} {s-2}}, $$ for $ 2 <| \nu | \leq s. $
Also taking into account that the derivative is represented as a linear combination of terms of the following form
$$ D ^ \alpha \exp \Big (\frac {- \| t \| ^ 2} {2} \Big) f (t) = \sum_ {0 \leq \beta \leq \alpha} D ^ { \alpha- \beta} \exp \Big (\frac {- \| t \| ^ 2} {2} \Big) D ^ \beta f (t), $$
and
$$ \Big | D ^ {\alpha- \beta} \exp \Big (\frac {- \| t \| ^ 2} {2} \Big) \Big | \leq c (\alpha- \beta, k) (1+ \| t \| ^ {| \alpha - \beta |}) \exp \Big (\frac {- \| t \| ^ 2} {2} \Big), $$
we get that
$$ \Big | D ^ \alpha \Big [\exp \Big (- \frac {\| t \| ^ 2} {2} + \sum_ {r = 1} ^ {s-3} \frac {\kappa_ {r + 2 } (it)} {(r + 2)!} \Big) - \exp \Big (\frac {- \| t \| ^ 2} {2} \Big) \sum_ {r = 1} ^ {s -3} P_r (it, \kappa_ \nu) \Big] \Big | $$
$$ \leq c (s, k) \rho_s (\| t \| ^ {s- | \alpha \|} + \| t \| ^ {3 (s-2) + | \alpha |}) \exp \Big (- \frac {1} {4} \| t \| ^ 2 \Big), $$
which completes the proof of the Lemma for the case when the covariance matrix $ D = I. $ In order to prove in the general case $ D \neq I $, it is necessary to consider the transformed sequence of random vectors
$ QZ_1, QZ_2,\dots,QZ_n, $ then $ \prod_{j = 1} ^ n \varphi_j (\theta_jQt) $ will be the characteristic function of the sum of random vectors $ Q (\sum_ {j = 1} ^ nZ_j ), $ and in this case,
the weighted sum of random vectors has a unit covariance matrix, and the Lemma holds for these random vectors.
\end {proof}
       \begin {lemma} \label {lemma4}
      Let $ X_1, X_2,\dots,X_n $ be independent random vectors with zero mean, unit covariance matrix and finite fourth absolute moment.
       Let $ l> 0 $ be such that for subset $ \mathbb {U} = \Big\{j: | \theta_j | \leq l/ \delta_j ^ 2 \Big\}, $ one has $ \sum \limits_ {j \in \mathbb {U}} \theta_j ^ 2 \geq 1/8. $
       If $ \delta_ \theta ^ 4 \leq (8k)^{-1}, $
then for $ \| t \| \leq (8 \sqrt {k} l)^{-1} $, it holds
       $$ \Big | D ^ \alpha \prod_ {j = 1} ^ n \varphi_j (\theta_j t) \Big | \leq c_1 (\alpha, k) (1+ \| t \| ^ {| \alpha |}) \exp \Big (- \frac {1} {48} \| t \| ^ 2 \Big). $$
       Where $ \delta_ \theta ^ 4 $ is defined in (\ref {d01}) and $ \varphi_j (t) $ in (\ref {df5}).
  \end {lemma}
  \begin {proof}
   The proof of Lemma follows the scheme of proofs of Lemma 2.2 \cite {KS} and Lemma 14.3 \cite {B}.
  For the random vector $ Z_j = Y_j- \mathbb {E} Y_j$, where $Y_j = X_j 1_X (\| \theta_j X_j \| \leq 1),$ ( see \ref{d2}, \ref{d3}), consider the expansion of the characteristic function of a Taylor series as $ t $ tends to zero,
  $$ \varphi_j (\theta_j t) = 1- \frac {1} {2} \theta_j ^ 2 \mathbb {E} \langle Z_j, t \rangle ^ 2 + \frac {1} {6} \theta_j ^ 3 \mathbb {E} \langle Z_j, t \rangle ^ 3 \xi, $$
  where $ | \xi | \leq 1. $
  For $ \| t \| \leq  (8 \sqrt {k} l)^{-1} $ and $ j \in \mathbb {U}$ one has
  $$ \theta_j ^ 2 \mathbb {E} \langle Z_j, t \rangle ^ 2 \leq \theta_j ^ 2 \| t \| ^ 2 \mathbb {E} \| Z_j \| ^ 2 \leq \frac {\mathbb {E} \| Z_j \| ^ 2 } {64k \delta_j ^ 4} \leq \frac {4\mathbb {E} \| X_j \| ^ 2} {64k \delta_j ^ 4}  = \frac {k} {16k\delta_j ^ { 4}} <1. $$
  Also
  $$ | \theta_j | ^ 3 \mathbb {E} | \langle Z_j, t \rangle | ^ 3 \leq | \theta_j | ^ 3 \| t \| ^ 3 \mathbb {E} \| Z_j \| ^ 3 \leq | \theta_j | ^ 2 \| t \| ^ 2 \frac {2 ^ 3 \mathbb {E} \| X_j \| ^ 3} {8 \sqrt {k} \delta_j ^ 2} $$
  $$ \leq | \theta_j | ^ 2 \| t \| ^ 2 \frac {2 ^ 3 \sqrt {k} \delta_j ^ 2} {8 \sqrt {k} \delta_j ^ 2} \leq | \theta_j | ^ 2 \| t \| ^ 2, $$
hence
  $$ 1- \frac {1} {2} \theta_j ^ 2 \mathbb {E} \langle Z_j, t \rangle ^ 2> 0. $$
  And
  $$ | \varphi_j (\theta_j t) | \leq 1- \frac {1} {2} \theta_j ^ 2 \mathbb {E} \langle Z_j, t \rangle ^ 2 + \frac {1} {6} | \theta_j | ^ 3 \mathbb {E} | \langle Z_j, t \rangle | ^ 3 $$
  $$ \leq \exp \Big (- \frac {1} {2} \theta_j ^ 2 \mathbb {E} \langle Z_j, t \rangle ^ 2 + \frac {1} {6} | \theta_j | ^ 3 \mathbb {E} | \langle Z_j, t \rangle | ^ 3 \Big) $$
  $$ \leq \exp \Big (- \frac {1} {2} \theta_j ^ 2 \mathbb {E} \langle Z_j, t \rangle ^ 2 + \frac {1} {6} \theta_j ^ 2 \| t \| ^ 2 \Big). $$
Further, note that
  $$ \exp \Big (\frac {1} {2} \theta_j ^ 2 \mathbb {E} \langle Z_j, t \rangle ^ 2- \frac {1} {6} | \theta_j | ^ 3 \mathbb {E} | \langle Z_j, t \rangle | ^ 3 \Big) $$
  $$ \leq \exp \Big (\frac {1} {2} (\mathbb {E} | \theta_j \langle Z_j, t \rangle | ^ 3) ^ {\frac {2} {3}} - \frac {1} {6} | \theta_j | ^ 3 \mathbb {E} | \langle Z_j, t \rangle | ^ 3 \Big) \leq \exp \Big (\frac {2} {3} \Big).
  $$
Now we denote the subset $ N_r = \{j_1, j_2 ,\dots, j_r \} $ as a subset of $ N = \{1,2,\dots,n \} $, consisting of $ r $ elements, for $ \| t \| \leq (8 \sqrt {k} l)^{-1} $ the following chain of inequalities holds
  $$
  \Big | \prod_ {j \in N \setminus N_r} \varphi_j (\theta_jt) \Big | \leq \Big | \prod_ {j \in \mathbb {U} \setminus N_r} \varphi_j (\theta_jt) \Big |
\leq \exp \Big (\sum_ {j \in \mathbb {U} \setminus N_r} \Big [- \frac {1} {2} \theta_j ^ 2 \mathbb {E} \langle Z_j, t \rangle ^ 2 + \frac {1} {6} \theta_j ^ 2 \| t \| ^ 2 \Big] \Big)
$$
$$
\exp \Big (\sum_ {j \in \mathbb {U} \cap N_r} \Big [- \frac {1} {2} \theta_j ^ 2 \mathbb {E} \langle Z_j, t \rangle ^ 2 + \frac {1} {6} \theta_j ^ 2 \| t \| ^ 2 \Big] - \sum_ {j \in \mathbb {U} \cap N_r} \Big [- \frac {1} {2 } \theta_j ^ 2 \mathbb {E} \langle Z_j, t \rangle ^ 2 + \frac {1} {6} \theta_j ^ 2 \| t \| ^ 2 \Big] \Big)
$$
$$
\leq \exp \Big (\sum_ {j \in \mathbb {U}} \Big [- \frac {1} {2} \theta_j ^ 2 \mathbb {E} \langle Z_j, t \rangle ^ 2 + \frac {1} {6} \theta_j ^ 2 \| t \| ^ 2 \Big] \exp \Big (\frac {2r} {3} \Big)
$$
  $$
  \leq \exp \Big (- \frac {1} {2} \sum_ {j = 1} ^ n \theta_j ^ 2 \mathbb {E} \langle Z_j, t \rangle ^ 2 + \frac {1} { 2} \sum_ {j \in N \setminus \mathbb {U}} \theta_j ^ 2 \mathbb {E} \langle Z_j, t \rangle ^ 2
  $$
  $$
  - \frac {1} {2} \sum_ {j \in N \setminus \mathbb {U}} \theta_j ^ 2 \| t \| ^ 2 + \frac {1} {2} \sum_ {j \in N \setminus \mathbb {U}} \theta_j ^ 2 \| t \| ^ 2 + \frac {1} {6} \| t \| ^ 2 \Big) \exp \Big (\frac {2r} {3} \Big)
  $$
  $$
  = \exp \Big (- \frac {1} {2} \langle Dt, t \rangle + \frac {1} {2} \Big [\langle D_2t, t \rangle- \| t \| ^ 2 \sum_ {j \in N \setminus \mathbb {U}} \theta_j ^ 2 \Big] + \frac {1} {2} \sum_ {j \in N \setminus \mathbb {U}} \theta_j ^ 2 \| t \| ^ 2 + \frac {1} {6} \| t \| ^ 2 \Big) \exp \Big (\frac {2r} {3} \Big).
  $$
Note that the norm of the covariance matrix satisfies the estimate $ \| D \|> 3/4, $ also, by definition, the inequality for the sum of the weight coefficients is satisfied $\sum_ {j \in N \setminus \mathbb {U}} \theta_j ^ 2 \leq 1/8. $ If we denote the matrix $ D_2 $ as the weighted sum of the covariance matrices of the random vectors $ X_j, j \in N \setminus \mathbb {U} $, then by Lemma \ref {lemma2} we obtain
  $$ \Big | \langle t, D_2t \rangle- \| t \| ^ 2 \sum_ {j \in N \setminus \mathbb {U}} \theta_j ^ 2 \Big | \leq 2k \sum_ {j \in N \setminus \mathbb {U}} \theta_j ^ 4 \delta_j ^ 4 \| t \| ^ 2 \leq \frac {1} {4} \| t \| ^ 2 , $$
  so
    $$
        \Big | \prod_ {j \in N \setminus N_r} \varphi_j (\theta_jt) \Big | \leq \exp \Big (\Big (- \frac {3} {8} + \frac {1} {8} + \frac {1} {16} + \frac {1} {6} \Big) \| t \| ^ 2 \Big) \exp \Big (\frac {2r} {3} \Big)
        $$
        \begin {equation} \label {pr}
        = \exp \Big (- \frac {1} {48} \| t \| ^ 2 \Big) \exp \Big (\frac {2r} {3} \Big).
        \nonumber
         \end {equation}
For $ \alpha = 0 $ the statement of the Lemma is proved.

Before proceeding to the proof of the case $ \alpha \neq 0 $, consider the modulus of the derivative of the characteristic function of the random vector $ Z_j $
\begin {equation} \label {der}
   \Big | D_m \varphi_j (\theta_j t) \Big | = | \theta_j | \Big | \mathbb {E} Z_ {j, m} \exp \Big (i \langle \theta_jt, Z_j \rangle \Big) \Big |,
\end {equation}
If a positive vector $ \beta $ satisfies the condition $ | \beta | = 1 $ then
$$ \Big | D _ {\beta} \varphi_j (\theta_j t) \Big | = | \theta_j | \Big | \mathbb {E} Z_ {j, \beta} \Big (\exp \Big (i \langle \theta_jt, Z_j \rangle \Big) -1 \Big) \Big | $$
$$ \leq | \theta_j | \mathbb {E} \Big | Z_ {j, \beta} \langle \theta_jt, Z_j \rangle \Big |
\leq \theta_j ^ 2 \| t \| \mathbb {E} \| Z_j \| ^ 2 \leq \theta_j ^ 2 \rho_2 (Z_j) \| t \|, $$
from (\ref {der}) we also get that for any vector $ \beta $ with $ | \beta | \geq 2 $ 
$$ \Big | D ^ \beta \varphi_j (\theta_j t) \Big | \leq | \theta_j | ^ {| \beta |} \mathbb {E} | Z_j ^ \beta | \leq 2 ^ {| \beta |} \theta_j ^ 2 \rho_2 (X_j), $$
finally we get that for any non-negative vector $ \beta> 0 $
\begin {equation} \label {pr2}
     \Big | D ^ \beta \varphi_j (\theta_j t) \Big | \leq c_2 (\alpha, k) \theta_j ^ 2 \rho_2 (X_j) \max \{1, \| t \| \}.\nonumber
\end {equation}
Now consider the positive vector $ \alpha> 0 $, according to the rule of differentiation of the product of functions, we obtain that
\begin {equation} \label {e6}
    D ^ \alpha \prod_ {j = 1} ^ n \varphi_j (\theta_j t) = \sum \prod_ {j \in N \setminus N_r} \varphi_j (\theta_j t) D ^ {\beta_1} \varphi_ { j_1} (\theta_ {j_1} t) \dots D ^ {\beta_r} \varphi_ {j_r} (\theta_ {j_r} t),
\end {equation}
where $ N_r = \{j_1,\dots,j_r \}, \; 1 \leq r \leq | \alpha |, $
$\beta_1, \beta_2, \beta_3 ,\dots, \beta_r $ vectors that meet the conditions
$ | \beta_j | \geq 1 \; \; (1 \leq j \leq r) $ and $ \sum_ {j = 1} ^ r \beta_j = \alpha. $
The number of multiplications in each of the $ n ^ \alpha $ terms of the expression (\ref {e6}) is
$$ \frac {\alpha_1! \dots \alpha_k!} {\prod_ {j = 1} ^ r \prod_ {i = 1} ^ k \beta_ {ji}! }, $$
where $ \alpha = (\alpha_1, \dots, \alpha_k) $ and $ \beta_j = (\beta_ {j_1}, \dots, \beta_ {j_k}), \; 1 \leq j \leq r. $ Each term in the expression (\ref {e6}) is bounded by the value
$$ \exp \Big (\frac {2r} {3} - \frac {1} {48} \| t \| \Big) \prod_ {j \in N_r} b_j, $$
where $ b_j = c_2 (\alpha, k) \rho_2 (X_j) \theta_j ^ 2 \max \{1, \| t \| \}, $
therefore from (\ref {e6}) we obtain
$$ \Big | D ^ \alpha \prod_ {j = 1} ^ n \varphi_j (\theta_j t) \Big | \leq \sum_ {1 \leq r \leq | \alpha |} c_3 (\alpha, r) \exp \Big (\frac {2r} {3} - \frac {1} {48} \| t \| \Big) \sum_r \prod_ {j \in N_r } b_j, $$
where the outer summation is over all $ r $ elements from $ N. $
It remains to evaluate the expression
$$ \sum_r \prod_ {j \in N_r} b_j \leq \Big (\sum_ {j = 1} ^ n b_j \Big) ^ r = (c_2 (\alpha, k) \rho_2 \max \{1, \| t \| \}) ^ r  =
(c_2 (\alpha, k) k) ^ r (1+ \| t \| ^ r), $$
which completes the proof.
  \end {proof}
 
  \begin {lemma} \label {lemma5}
         Let $ X_1, X_2,\dots,X_n $ be independent random vectors with zero mean, unit covariance matrix and finite fourth absolute moment. Let $ (\Theta_1, \Theta_2 ,\dots, \Theta_n) $ be a random vector uniformly distributed on the unit sphere $ S ^ {n-1}. $ Then with probability greater than $ 1- \grave {C_ {2}} ( \alpha, k) \exp \Big (- \grave {c} _ {2} (k) \frac {n} {\delta ^ 4} \Big) $ one has
$$ \int \limits _ {\frac {\beta (k) \sqrt n} {\delta ^ 2} \leq \| t \| \leq \frac {n} {\delta ^ 4}} \Big | D ^ \alpha \prod_ {j = 1} ^ n \varphi_j (\Theta_jt) \Big | dt \leq C_3 (\alpha, c, k) \frac {\delta ^ 4} {n}. $$
Where $ \varphi_j (t) $ is defined in (\ref {df5}) and $ \delta ^ 4 $ in (\ref {d01}).
  \end {lemma}
  \begin {proof}
   The proof of Lemma follows the scheme of the proof of Lemma 3.5 \cite {KS}.
  Let us estimate the modulus of the characteristic function of the truncated random vector 
  $$ 
  Y_j = X_j1_X (\| \Theta_jX_j \| \leq 1),
  $$
  ( see (\ref{d2}) ), denote $ g_j (t) $ as the characteristic function of the random vector $ X_j. $ First, using the Chebyshev inequality, we obtain the estimate
  $$
  \mathbb {E} 1_X (\| \theta_jX_j \| \geq 1) = P (\| \theta_jX_j \| \geq 1)  \leq \mathbb {E} \| X_j \| ^ 2 \theta_j ^ 2 = k \theta_j ^ 2,
  $$
then, applying a number of transformations, we obtain the following chain of inequalities for estimating the characteristic function
  $$
  \Big | \mathbb {E} \exp \Big (i \Theta_j \langle t, Y_j \rangle \Big) \Big | \Theta_j \Big | = \Big | \mathbb {E} \exp \Big (i \Theta_j \langle t, X_j \rangle 1_X (\| \Theta_jX_j \| \leq 1) \Big) \Big | \Theta_j \Big |
  $$
  $$
  = \Big | \mathbb {E} \exp \Big (i \Theta_j \langle t, X_j \rangle 1_X (\| \Theta_jX_j \| \leq 1) \Big) \Big (1_X (\| \Theta_jX_j \| \leq 1) + 1_X (\| \Theta_jX_j \|> 1) \Big) \Big | \Theta_j \Big |
  $$
  $$
    = \Big | \mathbb {E} \exp \Big (i \Theta_j \langle t, X_j \rangle \Big) 1_X (\| \Theta_jX_j \| \leq 1) + \mathbb {E} 1_X (\| \Theta_jX_j \|> 1) \Big | \Theta_j \Big |
  $$
  $$
  = \Big | \mathbb {E} \Big [\exp \Big (i \Theta_j \langle t, X_j \rangle \Big) - \exp \Big (i \Theta_j \langle t, X_j \rangle \Big) 1_X (\| \Theta_jX_j \|> 1) \Big] + \mathbb {E} 1_X (\| \Theta_jX_j \|> 1) \Big | \Theta_j \Big |
  $$
  $$
    = \Big | \mathbb {E} \Big [\exp \Big (i \Theta_j \langle t, X_j \rangle \Big) \Big] + \mathbb {E} 1_X (\| \Theta_jX_j \|> 1) \Big (1 - \exp \Big (i \Theta_j \langle t, X_j \rangle \Big) \Big | \Theta_j \Big |
  $$
  $$
  \leq | g_j (\Theta_jt) | + \Big | \mathbb {E} 1_X (\| \Theta_jX_j \|> 1) \Big (1- \exp \Big (i \Theta_j \langle t, X_j \rangle \Big) \Big | \Theta_j \Big |
  $$
    $$
  \leq | g_j (\Theta_jt) | + \mathbb {E} 1_X (\| \Theta_jX_j \|> 1) \Big | 1- \exp \Big (i \Theta_j \langle t, X_j \rangle \Big | \; \; \Big | \Theta_j 
  $$
  \begin {equation} \label {e7}
  \leq | g_j (\Theta_jt) | +2 \mathbb {E} 1_X (\| \Theta_jX_j \|> 1) \Big | \Theta_j \leq | g_j (\Theta_jt) | + 2k \Theta_j ^ 2.
  \end {equation}
  Next, we will show that for any $ r $ one has 
  $$ \mathbb {E} | g_j (\Theta_jr) | ^ 2 \leq 1 - c_1 \min \Big \{\frac {\| r \| ^ 2} {n},\frac{1}{ \delta_j ^ { 4}} \Big \}. $$
Let us denote $ X_j '$ as an independent copy of the random vector $ X_j $,
and define a random vector $ \hat {X} = X_j-X'_j $ with the corresponding distribution $ \hat {F} _x. $ Also denote $ f_n $ and $ J_n $ as the distribution density and characteristic function of the component of a random vector uniformly distributed on the unit sphere, $ \Theta_j. $ Further, changing the order of integration, we obtain
$$ \mathbb {E} | g_j (\Theta_jr) | ^ 2 = \int | g_j (rt) | ^ 2f_n (t) dt  = \int \Big [\int \exp \Big (it \langle t, y \rangle \Big) d \hat {F} _x \Big] f_n (t) dt
$$
$$ = \int \Big [\int f_n (t) \exp \Big (it \langle r, y \rangle \Big) dt \Big] d \hat {F} _x  = \int J_n (\langle r, x \rangle) d \hat {F} _x = \mathbb {E} J_n (\langle r, \hat {X} \rangle)
$$
$$ = \mathbb {E} J_n \Big (\| r \| \frac {\langle r, \hat {X} \rangle} {\| r \|} \Big). $$
Lemma 3.3 \cite {KS} implies that the estimate holds for the characteristic function of a random variable uniformly distributed on the unit sphere
$$ \mathbb {E} J_n \Big (\| r \| \frac {\langle r, \hat {X} \rangle} {\| r \|} \Big) \leq 1- \overline {c} _3 \mathbb {E} \min \Big \{\frac {\| r \| ^ 2} {n} \Big (\frac {\langle r, \hat {X} \rangle} {\| r \|} \Big) ^ 2,1 \Big \}.
$$
Let us define the random variable $ X '' $ as
$ X '' = \frac {\langle r, \hat {X} \rangle ^ 2} {2 \| r \| ^ 2} $ and $ \tau = \frac {\| r \| ^ 2} {n}, $
then
$$ \mathbb {E} X '' = \mathbb {E} \Bigg (\frac {\sum \limits_ {i = 1} ^ {k} r_i ^ 2 \hat {X} _i ^ 2 + 2 \sum \limits_ {1 \leq j <i \leq k} r_jr_i \hat {X} _j \hat {X} _i} {2 \| r \| ^ 2} \Bigg) = \frac {2 \| r \| ^ 2} {2 \| r \| ^ 2} = 1, $$
$$ \mathbb {E} (X '') ^ 2 = \frac {\langle r, \hat {X} \rangle ^ 4} {4 \| r \| ^ 4} \leq \mathbb {E} \frac {\| r \| ^ 4} {4 \| r \| ^ 4} \| \hat {X} \| ^ 4 \leq \frac {2 \delta_j ^ 4 + 6k ^ 2 } {4} \leq2 \delta_j ^ 4. $$
Let us show that $ \mathbb {E} \min \{\tau X '', 1 \} \geq \overline {c} _4 \min \{\tau, \delta_j ^ { -4} \}. $
Since the right-hand side increases with $ \tau $, it suffices to prove the inequality for $ \tau <(10 \delta_j ^ 4)^{-1}. $ From Lemma 3.1 \cite {KS} it follows
$\mathbb {E} 1_X (X '' \leq 10 \delta_j ^ 4) X '' \geq 4/5,$
so for $ 0 <\tau \leq (10 \delta_j ^ 4)^{-1} $ one has
$$
\mathbb {E} \min \{\tau X '', 1 \} \geq \mathbb {E} 1_X (X '' \leq 10 \delta_j ^ 4) \min \{\tau X '', 1 \}  = \tau \mathbb {E} 1_X (X '' \leq 10 \delta_j ^ 4) X ''> \frac {\tau} {2},
$$
hence
$$
\; \; \; \mathbb {E} | g_j (\Theta_jr) | ^ 2 \leq 1-c_1 \min \Big \{\frac {\| r \| ^ 2} {n},\frac{1}{ \delta_j ^ { 4}} \Big \}.
$$
Similarly, consider $ Z_j-Z'_j $, where $ Z'_j $ is an independent copy of $ Z_j $, also note that $ \mathbb {E} \Theta_j ^ 2 = n^{-1}$ 
and $\mathbb {E} \Theta_j ^ 4 \leq \hat{C}/n ^ {2}.$
And, using the inequality (\ref {e7}), we obtain the estimate
  $$ \mathbb {E} | \varphi_j (\Theta_jt) | ^ 2 = \mathbb {E} \Big [\mathbb {E} \exp \Big (i \Theta_j \langle t, Z_j-Z'_j \rangle \Big) \Big | \Theta_j \Big]
  = \mathbb {E} \Big [\mathbb {E} \exp \Big (i \Theta_j \langle t, Y_j-Y'_j \rangle \Big) \Big | \Theta_j \Big] 
  $$
  $$
  \leq \mathbb {E} \Big [\Big (| g_j (\Theta_jt) | + 2k \Theta_j ^ 2 \Big) \Big (| g_j (- \Theta_jt) | + 2k \Theta_j ^ 2 \Big ) \Big] 
  \leq \mathbb {E} | g_j (\Theta_j t) | ^ 2 + 2k \frac {1} {n} + 4k ^ 2\hat{C} \frac {1} {n ^ 2}
  $$
  $$ 
  \leq 1-c_1 \min \Big \{\frac {\| t \| ^ 2} {n},\frac{1}{ \delta_j ^ { 4}} \Big \} + 2k \frac {1} {n} + 4k ^ 2\hat{C} \frac {1} {n ^ 2}.
$$
Since $ \delta_j ^ 4 \geq k ^ 2 $ and in the region $ \| t \| ^ 2 \geq n/ \delta ^ 4 $, then the estimate holds
  $$ 
  \Big (1-c_1 \min \Big \{\frac {\| t \| ^ 2} {n},\frac{1}{ \delta_j ^ { 4}} \Big \} + 2k \frac {1} {n} + 4k\hat{C} \frac {1} {n ^ 2} \Big) ^ {- 1} \leq \Big(1- \frac {\overline{C}} {k ^ 2}\Big) ^ {- 1}. 
  $$
  In Lemma \ref {lemma4} it was shown that
  $$ \Big | D ^ \beta \varphi_j (\theta_j t) \Big | \leq c_2 (\alpha, k) \theta_j ^ 2 \rho_2 (X_j) \max \{1, \| t \| \}, $$
from this it immediately follows
  $$ \Big (\mathbb {E} \Big | D ^ \beta \varphi_j (\Theta_j t) \Big | ^ 2 \Big) ^ {\frac {1} {2}} \leq c_2 (\alpha, k) \rho_2 (X_j) \max \{1, \| t \| \} \Big (\mathbb {E} \Theta_j ^ 4 \Big) ^ {\frac {1} {2}} 
  $$
  $$\ = c_2 (\alpha, k) \rho_2 (X_j) \max \{1, \| t \| \} \frac {\sqrt {\hat{C}} } {n}.
  $$
Further, using Theorem 1 \cite {S}, we obtain
  $$
  \mathbb {E} \Big | \Big [\prod_ {j \in N_r} \varphi_j (\Theta_jt) \Big] D ^ {\beta_1} \varphi_ {j_1} (\Theta_ {j_1}) \dots D ^ {\beta_r} \varphi_ { j_r} (\Theta_ {j_r}) \Big |
  $$
  $$ \leq \Big [\prod_ {j \in N_r} \left (\mathbb {E} | \varphi_j (\Theta_jt) | ^ 2 \right) ^ {\frac {1} {2}} \Big] \left (\mathbb {E} | D ^ {\beta_1} \varphi_ {j_1} (\Theta_ {j_1}) | ^ 2 \right) ^ {\frac {1} {2}} \dots \left ( \mathbb {E} | D ^ {\beta_r} \varphi_ {j_r} (\Theta_ {j_r}) | ^ 2 \right) ^ {\frac {1} {2}}.
  $$
Let us denote the subset of indices $\mathfrak {G} = \Big\{j: \delta_j <2 \delta\Big \} 
$
and we come to the conclusion that
$$ 
\delta ^ 4 = \frac {1} {n} \sum \limits_ {j = 1} ^ {n} \delta_j ^ 4 \geq \frac {1} {n} \sum \limits_ {j \notin \mathfrak {G}} \delta_j ^ 4 \geq \frac {n- | \mathfrak {G} |} {n} 16 \delta ^ 4, 
$$
therefore,
$ | \mathfrak {G} | \geq  n/2. $ Due to this, the following chain of inequalities is valid 
$$ 
\mathbb {E} \Big | \prod \limits_ {j = 1} ^ {n} \varphi_j (\Theta_jt) \Big | \leq \prod \limits_ {j = 1} ^ {n} \Big (1-c_1 \min \Big \{\frac {\| t \| ^ 2} {n},\frac{1}{ \delta_j ^ { 4}} \Big \} + 2k \frac {1} {n} + 4k\hat{C} \frac {1} {n ^ 2} \Big) ^ {\frac {1} {2}}
$$
$$
\leq \prod \limits_ {j \in \mathfrak {G}} \Big (1-c_1 \min \Big \{\frac {\| t \| ^ 2} {n},\frac{1}{ \delta_j ^ { 4}} \Big \} + 2k \frac {1} {n} + 4k ^ 2\hat{C} \frac {1} {n ^ 2} \Big) ^ {\frac {1} {2}}
$$
$$
\leq \Big (1- \overline {c} _5 \min \Big \{\frac {\| t \| ^ 2} {n}, \frac{1}{ \delta ^ { 4}}\Big \} + 2k \frac {1} {n} + 4k ^ 2\hat{C} \frac {1} {n ^ 2} \Big) ^ {\frac {n} {4}}
$$
$$
\leq \exp \Big (\frac {k} {2} + \hat{C}k ^ 2 \Big) \exp \Big (- \overline {c} _6 \min \Big \{\| t \| ^ 2, \frac {n} {\delta ^ 4} \Big \} \Big). 
$$
Finally, we come to the conclusion
$$ \prod \limits_ {j \in N_r} \Big (\mathbb {E} | \varphi_j (\Theta_jt) | ^ 2 \Big) ^ {\frac {1} {2}} \leq \exp \Big (\frac {k} {2} + \hat{C}k ^ 2 \Big) \Big(1- \frac {\overline{C}} {k ^ 2}\Big) ^ {- r} \exp \Big (- \overline {c} _6 \min \Big \{\| t \| ^ 2, \frac {n} {\delta ^ 4} \Big \} \Big), $$
and similarly, as in the proof of Lemma \ref {lemma4}, we obtain the estimate
$$ 
\mathbb {E} \Big | D ^ \alpha \prod_j ^ n \varphi_j (\theta_j t) \Big | \leq c_3 (\alpha, k) (1+ \| t \| ^ {| \alpha |}) \exp \Big (- \overline {c} _6 \min \Big \{\| t \| ^ 2 , \frac {n} {\delta ^ 4} \Big \} \Big)
$$
and therefore
$$
\int \limits _ {\frac {\beta (k) \sqrt n} {\delta ^ 2} \leq \| t \| \leq \frac {n} {\delta ^ 4}} \mathbb {E} \Big | D ^ \alpha \prod_ {j = 1} ^ {n} \varphi_j (\Theta_jt) \Big | dt
$$
$$ \leq
\int \limits _ {\frac {\beta (k) \sqrt n} {\delta ^ 2} \leq \| t \| \leq \frac {n} {\delta ^ 4}} c_3 (\alpha, k) (1+ \| t \| ^ {| \alpha |}) \exp \Big (- \overline {c} _7 \min \Big \{\| t \| ^ 2, \frac {n} {\delta ^ 4} \Big \} \Big) dt
$$
$$
\leq \int \limits _ {\frac {\beta (k) \sqrt n} {\delta ^ 2} \leq \| t \|
\leq \frac {n} {\delta ^ 4}} c_3 (\alpha, k) (1+ \| t \| ^ {| \alpha |}) \exp \Big (- \overline {c} _8 \frac {n} {2 \delta ^ 4} \Big) dt
$$
$$ \leq \overline {C} _7 (\alpha, c, k) \exp \Big (- \overline {c} _9 (k) \frac {n} {\delta ^ 4} \Big).
$$
Using the Chebyshev inequality,
$$ P \Big (\int \limits _ {\frac {\beta (k) \sqrt n} {\delta ^ 2} \leq \| t \| \leq \frac {n} {\delta ^ 4}} \Big | D ^ \alpha \prod_ {j = 1} ^ {n} \varphi_j (\Theta_jt) \Big | dt \geq \sqrt {\overline {C} _7 (\alpha, k) \exp \Big ( - \overline {c} _9 (k) \frac {n} {\delta ^ 4} \Big)} \; \; \Big) $$
$$ \leq \sqrt {\overline {C} _7 (\alpha, k) \exp \Big (- \overline {c} _9 (k) \frac {n} {\delta ^ 4} \Big)}. $$
We see that there is a subset $ \mathbb {Q} _2 \subseteq S ^ {n-1} $ on the unit sphere with probability $$ P (\Theta \in \mathbb {Q} _2) \geq 1- \grave {C } _ {2} (\alpha, k) \exp \Big (- \grave {c} _ {2} (k) \frac {n} {\delta ^ 2} \Big), $$ such that for any vector of weight coefficients $ (\theta_1, \theta_2 ,\dots, \theta_n) \in \mathbb {Q} _2 $ the inequality
$$ \int \limits _ {\frac {\beta (k) \sqrt n} {\delta ^ 2} \leq \| t \| \leq \frac {n} {\delta ^ 4}} \Big | D ^ \alpha \prod_ {j = 1} ^ {n} \varphi_j (\theta_jt) \Big | dt \leq \sqrt {\overline {C} _7 (\alpha, k) \exp \Big (- \overline {c} _9 (k) \frac {n} {\delta ^ 4} \Big)}  \leq C_3 (\alpha, k) \frac {\delta ^ 4} {n}. $$
 
  \end {proof}
 
\begin {lemma} \label {lemma6}
      Let $ Z_1, Z_2,\dots,Z_n $ be independent random vectors with zero mean and finite absolute moment of order $ r + 2, \; r \geq1. $ If  the weighted covariance matrix $ D $ satisfies
      $ \| D \|> 3/4, $
then for any positive vector $ \alpha $ such that $ | \alpha | \leq 3r $ the following inequalities hold
$$ \int \Big | D ^ \alpha P_1 (it, \kappa_ \nu) \exp \Big (- \frac {1} {2} \langle Dt, t \rangle \Big) \Big | dt \leq C_6 (\alpha, k) \sum_ {| \nu | = 3} | \kappa_ \nu |,
$$
$$ \int \Big | D ^ \alpha P_r (it, \kappa_ \nu) \exp \Big (- \frac {1} {2} \langle Dt, t \rangle \Big) \Big | dt \leq C (\alpha, k, r) \rho_ {r + 2}, $$
for $ r \geq 2 $
, also
$$
\Big | D ^ \alpha \exp \Big (- \frac {1} {2} \langle Dt, t \rangle \Big) \Big | \leq C (\alpha, k) (1+ \| t \| ^ {| \alpha |}) \exp \Big (- \frac {3} {8} \| t \| ^ 2 \Big).
$$
Where the matrix $ D $ is defined in (\ref {d5}), $ \rho_ {r + 2} $ in (\ref {d7}) and the polynomial $ P_r (t, \kappa_ \nu) $ in (\ref{d13}), (\ref {d14}).
\end {lemma}
  \begin {proof}
   The proof of Lemma follows the scheme of the proof of Lemma 9.5 \cite {B}.
  Lemma 9.3 \cite {B} implies that
  $$ \Big | D ^ \alpha \exp \Big (- \frac {1} {2} \langle Dt, t \rangle \Big) \Big | \leq C (\alpha, k) (1+ \| t \| ^ {| \alpha |}) \exp \Big (- \frac {3} {8} \| t \| ^ 2 \Big), $$
  now we show that for any positive vector
  $ \alpha $ such that $ 0 \leq | \alpha | \leq 3r $ one has
  $$ \Big | D ^ \alpha P_r (z, \kappa_ \nu) \Big | \leq C_3 (\alpha, r) (1+ \rho_2 ^ {r-1}) (1+ \| z \| ^ {3r- | \alpha |}) \rho_ {r + 2} .$$
Differentiating the polynomial of the asymptotic expansion, we obtain the following representation
  $$ D ^ \alpha P_r (z, \kappa_ \nu) = \sum_ {m = 1} ^ r \frac {1} {m!} \sum_ {i_1,\dots,i_m: \sum_ {j = 1} ^ mi_j = r} \Bigg [\sum _ {\nu_1 ,\dots, \nu_m: | \nu_j | = i_j + 2} \frac {\kappa _ {\nu_1} \dots \kappa _ {\nu_m}} {\nu_1!. \dots \nu_m!} \Bigg]
$$
$$
\times \frac {(\nu_1 + \dots + \nu_m)!} {((\nu_1 + \dots + \nu_m- \alpha)!)} z ^ {\nu_1 + \dots + \nu_m- \alpha},
$$
where  $\kappa_ \nu=\sum_ {j = 1} ^ n \theta_j ^ {| \nu |} \kappa_ \nu (Z_j),$ ( see (\ref{d11}), (\ref{d13}), (\ref{d14})  ).
Estimating the cumulants in this expression through the corresponding absolute moments( see (\ref{d10}) ), we come to the conclusion
$$
| \kappa _ {\nu_1} \dots \kappa _ {\nu_m} | \leq c_1 (\nu_1) \rho_ {i_1 + 2} \dots c_1 (\nu_m) \rho_ {i_m + 2}
$$
$$
\leq c_1 (\nu_1) \dots c_1 (\nu_m) \rho_2 ^ {\frac {i_1 + 2 + \dots + i_m + 2} {2}} \frac {\rho_ {i_1 + 2} \dots \rho_ {i_m + 2}} {\rho_2 ^ {\frac {i_1 + 2} {2}} \dots \rho_2 ^ {\frac {i_m + 2} {2}}}
$$
$$
\leq c_1 (\nu_1) \dots c_1 (\nu_m) \rho_2 ^ {\frac {r} {2} + m} \Bigg (\frac {\rho_ {r + 2}} {\rho_2 ^ {\frac {r + 2} {2}}} \Bigg) ^ {\frac {i_1 + \dots + i_m} {r}} = c_1 (\nu_1) \dots c_1 (\nu_m) \rho_2 ^ {m-1} \rho_ {r + 2},
$$
using the well-known inequality $ \| t \| ^ a \leq \| t \| ^ b + \| t \| ^ c $ for $ 0 \leq b \leq a \leq c $, we come to the fact that
$$ | D ^ \alpha P_r (z, \kappa_ \nu) | \leq C_3 (\alpha, r)(1+ \rho_2 ^ {r-1}) (1+ \| z \| ^ {3r- | \alpha |}) \rho_ {r + 2}. $$
Further, since
$$ D ^ \alpha \exp \Big (- \frac {1} {2} \langle Dt, t \rangle \Big) P_r (it, \kappa_ \nu) = \sum_ {0 \leq \beta \leq \alpha} D ^ {\alpha- \beta} \exp \Big (- \frac {1} {2} \langle Dt, t \rangle \Big) D ^ \beta P_r (it, \kappa_ \nu), $$
we obtain that for any $ r \geq 1 $
$$
\int \Big | D ^ \alpha P_r (it, \kappa_ \nu) \exp \Big (- \frac {1} {2} \langle Dt, t \rangle \Big) \Big | dt \leq C (\alpha, k, r)  (1+ \rho_2 ^ {r-1}) \rho_ {r + 2}. 
$$
Separately, it is necessary to consider the case when $ r = 1 $
$$
\int \Big | D ^ \alpha P_1 (it, \kappa_ \nu) \exp \Big (- \frac {1} {2} \langle Dt, t \rangle \Big) \Big | dt 
$$
$$ 
= \int \Big | \sum_ {0 \leq \beta \leq \alpha} D ^ {\alpha- \beta} P_1 (it, \kappa_ \nu) D ^ {\beta} \exp \Big ( - \frac {1} {2} \langle Dt, t \rangle \Big) \Big | dt
$$
$$
= \int \Big | \sum_ {0 \leq \beta \leq \alpha} D ^ {\alpha- \beta} \sum_ {| \nu | = 3} \frac {\kappa _ {\nu}} {\nu!} (it) ^ {\nu} D ^ {\beta} \exp \Big (- \frac {1} {2} \langle Dt, t \rangle \Big) \Big | dt
$$
$$
\leq \int \Big | \sum_ {0 \leq \beta \leq \alpha} D ^ {\alpha- \beta} \sum_ {| \nu | = 3} \frac {\kappa _ {\nu}} {\nu!} (it) ^ {\nu} \Big | C (\beta, k) (1+ \| t \| ^ {| \beta | }) \exp \Big (- \frac {3} {8} \| t \| ^ 2 \Big)) dt
$$
$$
\leq \sum_ {| \nu | = 3} | \kappa_ \nu | \Big [\int \Big | \sum_ {0 \leq \beta \leq \alpha} D ^ {\alpha- \beta} \frac {(it) ^ {\nu}} {\nu!} \Big | C (\beta, k) ( 1+ \| t \| ^ {| \beta |}) \exp \Big (- \frac {3} {8} \| t \| ^ 2 \Big) dt \Big]
$$
$$ \leq \sum_ {| \nu | = 3} C_6 (\alpha, \nu, k) | \kappa_ \nu | \leq C_6 (\alpha, k) \sum_ {| \nu | = 3} | \kappa_ \nu |.
$$
  \end {proof}
  \begin {lemma} \label {lemma7}
Let $ \Theta $ be a random vector uniformly distributed on the unit sphere $ S ^ {n-1} $, $ \nu $ is a positive vector with $ | \nu | = 3 $, then for any $ t> 0 $,
$$ P \Big (\Big | \sum \limits_ {j = 1} ^ n \Theta_j ^ 3 \mu _ {\nu} (Z_j) \Big | \geq t \frac {\delta ^ 4} {n} \Big) \leq \grave {C} _3 \exp \Big (- \grave {c} _3t ^ \frac {2} {3} \Big) $$
and
$$ P \Big (\sum \limits_ {j = 1} ^ {n} \delta_j ^ 4 \Theta_j ^ 4 \geq t \frac {\delta ^ 4} {n} \Big) \leq \grave {C } _4 \exp \Big (- \grave {c} _4 \sqrt t \Big). $$
Where $ \delta_j ^ 4 $, $ \delta ^ 4 $, $ \mu_ \nu (Z_j) $ are defined in (\ref {d02}), (\ref {d01}), (\ref {d9.1} ), respectively.
  \end {lemma}
  \begin {proof}
 
The proof follows the scheme of the proof of Lemma 4.1 \cite {KS}. Let $ \Gamma_1, \Gamma_2, ,\dots, \Gamma_n $ be a sequence of independent random variables with standard normal distribution and $ Z $ as a random variable independent of $ \Gamma_j $ for $ j = 1,\dots,n $, which has a chi-square with $ n $ degrees of freedom.
  Using the properties of the normal distribution, one can show that the representation $ \Theta_j = \Gamma_j/\sqrt Z $ is holds for $ j = 1,\dots,n. $ Applying a number of transformations, we obtain
 $$ P \Big (\Big | \sum \limits_{j = 1}^{n} \mu_{\nu} (Z_j) \Theta_j^3 \Big | \geq t \frac{\delta^4}{n} \Big) = P \Big (\Big | \sum \limits_{j = 1}^{n} \frac{\mu_{\nu} (Z_j) \Gamma_j^3}{Z^{\frac{3}{2}}} \Big | \geq t \frac{\delta^4}{n} \Big)
 $$
  $$ 
  \leq P \Big (\Big | \sum \limits_{j = 1}^{n} \mu_{\nu} (Z_j) \Gamma_j^3 \Big | \geq \frac{t \sqrt n \delta^2}{4} \Big) + \overline{C}_1 \exp \Big (- \overline{c}_1n \Big).
  $$
  Due to the fact that $\mathbb {E} \exp (c_2 \Gamma_j ^ 2) \leq2,$ where $ c_2 $ is an absolute constant, therefore the random variable $
  Y = \sum_ {j = 1} ^ {n} \mu _ {\nu} (Z_j) \Gamma_j ^ 3 
  $
  belongs to class $ \psi_ {2/3} $ see Sect. 2 and 3 in \cite{A}. And, therefore, we can apply the inequality for the moments Sect 3 in \cite {A}, which asserts that for any $ p \geq2 $ the following estimate holds
$$
\Big (\mathbb {E} | Y | ^ p \Big) ^ {\frac {1} {p}} \leq \overline{C}_2 p ^ {\frac {3} {2 }} \Big (\sum \limits_ {j = 1} ^ {n} \mu _ {\nu} (Z_j) ^ 2 \Big) ^ {\frac {1} {2}}.
$$
Using the inequality for $ \nu = \nu_1 + \nu_2 $ with $ | \nu_1 | = 1, \; | \nu_2 | = 2 $
 $$
(\mu _ {\nu} (Z_j)) ^ 2 = | \mathbb {E} Z_j ^ {\nu_1} Z_j ^ {\nu_2} | ^ 2 \leq (\mathbb {E} Z_j ^ {2 \nu_1} ) (\mathbb {E} Z_j ^ {2 \nu_2})  \leq2 ^ 2 \mathbb {E} X_j ^ {2 \nu_1} 2 ^ 4 \rho_4 (X_j) = 2 ^ 6 \delta_j ^ 4,
  $$
  we get that
  $$\overline{C}_{2}p^{\frac{3}{2}} \Big (\sum \limits_{j = 1}^{n} (\mu_{\nu} (Z_j))^2 \Big)^{\frac{1}{2}}\leq
\overline{C}_{3}p^{\frac{3}{2}} \Big (\sum \limits_{j = 1}^{n} \delta_j^4 \Big)^{\frac{1}{2} } =
\overline{C}_{3}p^{\frac{3}{2}} \sqrt n \delta^2.
$$
Applying the Chebyshev inequality, we come to the conclusion that
$$ 
P \Big (\Big | \sum \limits_{j = 1}^{n} \mu_{\nu} (X_j) \Gamma_j^3 \Big | \geq t \sqrt n \delta^2 \Big) \leq \frac{\mathbb{E} | Y |^p}{(t \sqrt n \delta^2)^p} \leq \Big (\frac{\overline{C}_{3} p^{\frac{3}{2}} }{t} \Big)^p.
$$
Put $ p = \Big(\frac{1}{2}\overline{C}_{3}^{-1}t\Big)^{\frac{2}{3}}, $ if we consider $  t \geq 2^{\frac{5}{2}} \overline{C}_{3} , $ then we get that $ p \geq 2. $ Note that
$$
\Big (\frac{\overline{C}_{3} p^{\frac{3}{2}} }{t} \Big)^p=\exp\Big(-\Big(\frac{1}{2}\overline{C}_{3}^{-1}t\Big)^{\frac{2}{3}}\ln(2)\Big),
$$ 
therefore we get that for any  $ t \geq \overline{c}_3 $
$$
P \Big (\Big | \sum \limits_{j = 1}^{n} \mu_{\nu} (X_j) \Theta_j^3 \Big | \geq t \frac{\delta^4}{n} \Big) \leq \exp \Big (- \overline{c}_3t^{\frac{2}{3}} \Big) + \overline{C}_1 \exp (- \overline{c}_1n),
$$
as
$$ 
\Big | \sum \limits_{j = 1}^{n} \mu_{\nu} (X_j) \Theta_j^3 \Big | \leq
\Big (\sum \limits_{j = 1}^{n} (\mu_{\nu} (X_j))^2 \Theta_j^4 \Big)^{\frac{1}{2}}\leq \Big (\sum \limits_{j = 1}^{n} 2^6\delta_j^4 \Theta_j^4 \Big)^{\frac{1}{2}} \leq 2^3\sqrt n \delta^2,$$
therefore, for any $ t \geq 0 $ one has
$$ P \Big (\Big | \sum \limits_{j = 1}^{n} \beta_j \Theta_j^3 \Big | \geq t \frac{\delta^4}{n} \Big) \leq \breve {C}_3 \exp \Big (- \breve {c}_3 t^\frac{2}{3} \Big). $$
Similarly for the random variable
$ Y_2 = \sum_ {j = 1} ^ {n} \delta_j ^ 4 (\Gamma_j ^ 4-3) $ we get
$$ P \Big (\sum \limits_{j = 1}^{n} \delta_j^4 \Theta_j^4 \geq 12 \frac{\delta^4}{n} + t \frac{\delta^4 }{n} \Big) \leq P \Big (Y_2 \geq t\frac{n \delta^4}{4} \Big) + \overline{C}_{4} \exp (- \overline{c}_{4}n) $$ 
and reusing the moment inequality (section 3 in \cite {A}) $$ \Big (\mathbb{E} | Y_2 |^p \Big)^{\frac{1}{p}} \leq \overline{C}_{5} p^2 n \delta^4. $$
Similarly, we come to the conclusion that
$$ P \Big (\sum \limits_{j = 1}^{n} \delta_j^4 \Theta_j^4 \geq t \frac{\delta^4}{n} \Big) \leq \overline{C }_{6} \exp \Big (- \overline{c}_{5} \sqrt t \Big) + \overline{C}_{4} \exp (- \overline{c}_{4}n),
t> 15. $$ 
Also from the fact that $ \sum_ {j = 1} ^ {n} \delta_j ^ 4 \Theta_j ^ 4 \leq n \delta ^ 4, $ it follows that for any $ t> 0 $
$$ P \Big (\sum \limits_ {j = 1} ^ {n} \delta_j ^ 4 \Theta_j ^ 4 \geq t \frac {\delta ^ 4} {n} \Big) \leq \breve {C } _4 \exp (- \breve {c} _4 \sqrt t). $$
  \end {proof}
          \section {Proof of the main theorem} \label {3}
        \begin {proof}
         Assume that $ \delta _ {\theta} ^ 4 \leq (8k)^{-1}.$
We split the original integral into several terms, for this we add and subtract the term with the distribution of the sum of truncated random variables defined in (\ref {df1}), (\ref {df2}), then estimate each term separately
$$ \Big | \int \limits_ {B} d (F_X- \Phi) \Big | \leq \Big | \int \limits_ {B} d (F_X-F_Y) \Big | + \Big | \int \limits_ {B} d (F_Y- \Phi) \Big | .$$
To estimate the first term in the sum, we use the following transformation
$$ \Big | \int \limits_ {B} d (F_X-F_Y) \Big | \leq \Big | \int \limits_ {B} d \Big (\sum_ {j = 1} ^ nF_ {X_j} \ast\dots\ast F_ {X_ {j-1}} \ast (F_ {X_ {j}} - F_ {Y_ {j}}) \ast F_ {Y_ {j + 1}} \ast\dots\ast F_ {Y_ {n}} \Big) \Big | $$
$$ \leq \sum_ {j = 1} ^ n \Big | \int \limits_ {B} d (F_ {X_j} -F_ {Y_j}) \Big | \leq \sum_ {j = 1} ^ nP (\| \theta_jX_j \|> 1) \leq \sum_ {j = 1} ^ n \delta_j ^ 4 \theta_j ^ 4 = \delta_ \theta ^ 4, $$
where $ F_ {X_j}, F_ {Y_j} - $ are defined (\ref {df4}), (\ref {df4_1}) and $ "\ast" $ denotes the function convolution.
To estimate the second term, we additionally split this integral into the sum of three integrals, adding and subtracting new terms
$$ \Big | \int \limits_ {B} d (F_Y- \Phi) \Big | = \Big | \int \limits_ {B_n = B + \{A_n\}} d (F_Z- \Phi_ {A_n, I}) \Big | $$
$$ \leq \Big | \int \limits_ {B_n} d (F_Z- \Phi_ {0, D}) \Big | + \Big | \int \limits_ {B_n} d (\Phi- \Phi_ {0, D}) \Big | + \Big | \int \limits_ {B_n} d (\Phi_ {A_n, I} - \Phi) \Big |. $$
To estimate the last integral, we show that due to $\mathbb {E}X_j=0,\;j=1,\dots,n,$ the following inequality holds:
$$ \| \mathbb {E} \theta_jY_ {j} \| ^ 2 = \sum_ {i = 1} ^ k \Big (\mathbb {E} \theta_j (Y_ {ji} -X_ {ji}) \Big) ^ 2 \leq
 \sum_ {i = 1} ^ k \Big (\mathbb {E} \| \theta_jX_ {j} \| 1_X (\| \theta_jX_ {j} \|> 1) \Big) ^ 2 $$
$$ \leq k \Big (\mathbb {E} \| \theta_jX_ {j} \| ^ 4 \Big) ^ 2 = k (\theta_j ^ 4 \delta_j ^ 4) ^ 2 .$$
From this it follows that the norm of the weighted mathematical expectation is bounded by the value
$$ \| A_n \| \leq \sum_ {j = 1} ^ n \| \mathbb {E} \theta_jY_j \| \leq \sqrt {k} \sum_ {j = 1} ^ n \theta_j ^ 4 \delta_j ^ 4 <\frac {1} {8 \sqrt {k}}. $$
From Theorem 4 \cite {U} we obtain that
$$ \Big | \int \limits_ {B_n} d (\Phi_ {A_n, I} - \Phi) \Big | \leq \sqrt {\frac {2} {\pi}} \| A_n \| \leq C_1 \sqrt k \delta_ \theta ^ 4.$$
To estimate the next integral from the sum, we use Theorem 3 \cite {U} and Lemma \ref {lemma2}
$$ \Big | \int \limits_ {B_n} d (\Phi- \Phi_ {0, D}) \Big | \leq C_2 \| I-D \| _2 = C_2 \Big (\sum_ {i, l} (d_ {il} -v_ {il}) ^ 2 \Big) ^ {\frac {1} {2}} $$
$$ \leq C_2 \Big (\sum_ {i, l} (2 \delta_ \theta ^ 4) ^ 2 \Big) ^ {\frac {1} {2}} \leq 2kC_2 \delta_ \theta ^ 4 .$$
It remains to estimate the last integral, for this we use the most important inequality from Corollary 11.5 [2]
$$ 
\Big | \int \limits_ {B_n} d (F_Z- \Phi_ {0, D}) \Big | \leq
C_3 \int \Big | d ((F_Z- \Phi_ {0, D}) \times K_ \epsilon) \Big | + \hat {C} _1 (k) \epsilon, 
$$
where $ \epsilon = 4a \sqrt k \frac {\delta ^ 4} {n} $ and $ a = 2 ^ {- \frac {1} {3}} k ^ {\frac {5} {6}} $, and $ K_ \epsilon (x) $ is a kernel function (more details 13.8-13.13 in \cite {B}),
the most important property of this function is that its characteristic function $ \widehat {K} (t) = 0 $ is equal to zero for $ \| t \|> n/ \delta ^ 4. $
By Lemma 11.6 \cite {B}, from estimating the difference of distributions, we can go to estimating the difference of the corresponding characteristic functions
$$ \int \limits \Big | d (F_Z- \Phi_ {0, D}) \times K_ \epsilon) \Big | \leq \hat {C} _2 (k) \max \limits_ {0 \leq | \alpha + \beta | \leq k + 1} \int \Big | D ^ \alpha (\hat {F} _Z- \hat {\Phi} _ {0, D}) (t) D ^ \beta \hat {K} _ \epsilon (t) \Big | dt.$$
Since $| D ^ \beta \hat {K} _ \epsilon (t) | \leq \hat {c}, $ we get that
$$ \int \limits \Big | D ^ \alpha (\hat {F} _Z- \hat {\Phi} _ {0, D}) (t) D ^ \beta \hat {K} _ \epsilon (t) \Big | dt \leq \hat {C} _3 (k) \int \limits _ {\| t \| \leq \frac {n} {\delta ^ 4}} \Big | D ^ \alpha (\hat {F} _Z- \hat {\Phi} _ {0, D}) (t) \Big | dt .$$
Denote $ E_n = c_1 (k, k + 3) \min \Big\{\eta_ {k + 3} ^ {- \frac {1} {k + 3}}, \eta_ {k + 3} ^ {- \frac {1} {k + 1}} \Big\}, $ where $   \eta_ {k + 3} = \sum_ {j = 1} ^ n \rho_{k+3} (Q \theta_jZ_j) $ ( see (\ref {d8}) ) and $ c_1 (k, k + 3) $ is a constant from the statement of Lemma \ref {lemma3}, further, we add and subtract terms of the asymptotic expansion of the logarithm of the characteristic function. Considering that, by definition, $ \hat {F} _Z = \prod_ {j = 1} ^ n \varphi_j (\theta_jt)$ ( see (\ref {df5}) ), we can split the integral into the sum of several terms by dividing the region of integration into several parts, similarly to inequality (71) \cite {Saz}
$$ 
\int \limits _ {\| t \| \leq \frac {n} {\delta ^ 4}} \Big | D ^ \alpha (\hat {F} _Z- \hat {\Phi} _ {0, D}) (t) \Big | dt
$$
$$ 
\leq
\int \limits _ {\| t \| \leq \sqrt {\frac {4} {5}} E_n} \Big | D ^ \alpha \Big [\prod_ {j = 1} ^ n \varphi_j (\theta_jt) - \exp \Big (- \frac {1} {2} \langle Dt, t \rangle \Big) \sum_ { r = 0} ^ {k} P_r (it, \kappa_v) \Big] \Big | dt
$$
$$ + \int \limits _ {\sqrt {\frac {4} {5}} E_n \leq \| t \| \leq \frac {\sqrt n} {1600 \sqrt k \delta ^ 2}} \Big | D ^ \alpha \prod_ {j = 1} ^ n \varphi_j (\theta_jt) \Big | dt + \int \limits _ {\frac {\sqrt n} {1600 \sqrt k \delta ^ 2} \leq \| t \| \leq \frac {n} {\delta ^ 4}} \Big | D ^ \alpha \prod_ {j = 1} ^ n \varphi_j (\theta_jt) \Big | dt
$$
$$ + \int \limits _ {\sqrt {\frac {4} {5}} E_n \leq \| t \|} \Big | D ^ \alpha \exp \Big (- \frac {1} {2} \langle Dt, t \rangle \Big) \Big | dt
+ \int \Big | D ^ \alpha \sum_ {r = 1} ^ {k} P_r (it, \kappa_ \nu) \exp \Big (- \frac {1} {2} \langle Dt, t \rangle \Big) \Big | dt 
$$
$$ = I_1 + I_2 + I_3 + I_4 + I_5 .$$
Then, by Lemma \ref {lemma2}, we obtain
$$ \det Q \leq \| Q ^ 2 \| ^ {\frac {k} {2}} \leq \Big (\frac {4} {3} \Big) ^ {\frac {k} {2 }}, $$
$$
\| Qt \| \geq \| D \| ^ {- \frac {1} {2}} \| t \| \geq \Big (\frac {4} {5} \Big) ^ {\frac {1} {2}} \| t \|,
$$
$$
\Big\{\| Qt \| \leq \sqrt {\frac {4} {5}} E_n \Big\} \subset \Big\{\| t \| \leq c_1 (k, k + 3) \min \Big\{\eta_ {k + 3} ^ {- \frac {1} {k + 3}}, \eta_ {k + 3} ^ {- \frac {1 } {k + 1}} \Big\} \Big\},
$$
where the matrix $ Q $ is defined in (\ref {d6}). Also, taking into account that for any $ s \geq4 $ from (\ref {d8b}) it follows
$$
\eta_s \leq 2 ^ s \| Q \| ^ s \sum_ {j = 1} ^ n \rho_4 (\theta_j X_j) = \| Q \| ^ s2 ^ s \delta_ \theta ^ 4 \leq 2 ^ s \Big (\frac {4} {3} \Big) ^ {\frac {s} {2}} \delta_ \theta ^ 4.
$$
Substituting $ t = Qt $, $ s = k + 3 $ and using Lemma \ref {lemma3}, we come to the fact that
$$ 
I_1 = \int \limits _ {\| Qt \| \leq \sqrt {\frac {4} {5}} E_n} \Big | D ^ \alpha \Big [\prod_ {j = 1} ^ n \varphi_j (\theta_jQt) - \exp \Big (\frac { - \| t \| ^ 2} {2} \Big) \sum_ {r = 0} ^ {k} P_r (iQt, \kappa_v) \Big] \det Q \Big | dt
$$
$$
\leq \int \limits _ {\| t \| \leq E_n} \Big (\frac {4} {3} \Big) ^ {\frac {k} {2}} \Big | D ^ \alpha \Big [\prod_ {j = 1} ^ n \varphi_j (\theta_jQt) - \exp \Big (\frac {- \| t \| ^ 2} {2} \Big) \sum_ {r = 0} ^ {k} P_r (iQt, \kappa_v) \Big] \Big | dt
$$
$$
\leq \int \Big (\frac {4} {3} \Big) ^ {\frac {k} {2}} c_2 (k, k + 3) \eta_ {k + 3} (\| t \| ^ {k + 3- | \alpha |} + \| t \| ^ {3 (k + 1) + | \alpha |}) \exp \Big (- \frac {1} {4} \| t \| ^ 2 \Big) dt
$$
$$
\leq \hat {C_1} (\alpha, k) \eta_ {k + 3} \leq 2 ^ {k + 3} \Big (\frac {4} {3} \Big) ^ {\frac {k + 3} {2}} \hat {C_1} (\alpha, k) \delta_ \theta ^ 4 \leq C_1 (\alpha, k) \delta_ \theta ^ 4.
$$
Then, we denote the subset $ 
\mathfrak {G} = \Big\{j: \; \delta_j ^ 2 <5 \delta ^ 2 \Big\},
$
the following holds
$$
\delta ^ 4 = \frac {1} {n} \sum \limits_ {j = 1} ^ {n} \delta_j ^ 4 \geq
\frac {1} {n} \sum \limits_ {j \notin \mathfrak {G}} \delta_j ^ 4> \frac {n- | \mathfrak {G} |} {n} (5 \delta ^ 2) ^ 2,
$$
hence $| \mathfrak {G} |/n> 24/25> 4/5.$
Using Lemma 3.2 \cite {KS}, we obtain that with probability greater than $ 1- \grave {C} _1 \exp (- \grave {c} _1n) $ a random vector $ \Theta $ uniformly distributed on the unit sphere $ S ^ {n -1} $ satisfies the following condition $$
\sum \limits_ {j \in \mathbb {U}} \Theta_j ^ 2> \frac {1} {8}, \mathbb {U} = \Big\{j \in \mathfrak {G}: | \Theta_j | <\frac {40} {\sqrt n}\Big \}.
$$
If we set $ l = 200 \delta ^ 2/\sqrt n,$ then one has $40 /\sqrt n \leq l / \delta_j ^ 2, $ for $ j \in \mathfrak {G} $. Therefore, using Lemma \ref {lemma4}, we obtain that there is a subset of $ \mathbb {Q} _1 $ with 
$ \lambda_{n-1} (\mathbb {Q} _1) \geq 1- \grave {C} _1 \exp (- \grave {c} _1n) $
such that 
$$
I_2 = \int \limits _ {\sqrt \frac {4} {5} E_n \leq \| t \| \leq \frac {\sqrt n} {1600 \sqrt {k} \delta ^ 2}} \Big | D ^ \alpha \prod_ {j = 1} ^ n \varphi_j (\theta_jt) \Big | dt 
$$
$$
\leq \int \limits _ {\sqrt \frac {4} {5} E_n \leq \| t \| \leq \frac {\sqrt n} {1600 \sqrt {k} \delta ^ 2}} c_1 (\alpha, k) (1+ \| t \| ^ {| \alpha |}) \exp \Big ( - \frac {1} {48} \| t \| ^ 2 \Big) dt
$$
$$
\leq \int \limits \Big (\| t \| ^ {- 1} \Big (\frac {4} {5} \Big) ^ {\frac {1} {2}} c_1 (k, k + 3) \min \Big\{\eta_ {k + 3} ^ {- \frac {1} {k + 3}}, \eta_ {k + 3} ^ {- \frac {1} {k + 1}} \Big\} \Big) ^ {- k-3}
$$
$$\times c_1 (\alpha, k) (1+ \| t \| ^ {| \alpha |}) \exp \Big (- \frac {1} {48} \| t \| ^ 2 \Big) dt
$$

$$
\leq \hat {C} _2 (\alpha, k) \min \Big\{\eta_ {k + 3} ^ {- \frac {1} {k + 3}}, \eta_ {k + 3} ^ {- \frac {1} {k + 1}} \Big\} ^ {- k-3}
= \hat {C} _2 (\alpha, k) \eta_ {k + 3} \max \Big\{1, \eta_ {k + 3} ^ {\frac {2} {k+1}} \Big\}
$$
$$
\leq \hat {C} _2 (\alpha, k) 2 ^ {k + 3} \Big (\frac {4} {3} \Big) ^ {\frac {k + 3} {2}} \delta_ \theta ^ 4 \Big(1+ \eta_ {k + 3} ^ {\frac {2} {k + 1}}\Big) \leq C_2 (\alpha, k) \delta_ \theta ^ 4.
$$
By Lemma \ref {lemma5} there exists a subset $ \mathbb {Q} _2 $ with  $\lambda_{n-1} (\mathbb {Q} _2) \geq 1- \grave {C} _ {2} (\alpha,k) \exp \Big (- \grave {c} _ {2} (k) \frac {n} {\delta ^ 2} \Big) $ such that for any vector of weight coefficients $ (\theta_1, \theta_2 ,\dots, \theta_n ) \in \mathbb {Q} _2 $ one has
$$
I_3 = \int \limits_ {\frac {\sqrt n} {1600 \sqrt {k} \delta ^ 2} \leq \| t \| \leq \frac {n} {\delta ^ 4}} \Big | D ^ \alpha \Big [\prod_ {j = 1} ^ n \varphi_j (\theta_jt) \Big] \Big | dt \leq C_3 (\alpha, k) \frac {\delta ^ 4} {n}.
$$
By Lemma \ref {lemma6}
$$
I_4 = \int \limits _ {\sqrt {\frac {4} {5}} E_n \leq \| t \|} \Big | D ^ \alpha \exp \Big (- \frac {1} {2} \langle Dt, t \rangle \Big) \Big | dt
$$
$$
\leq \int \limits _ {\sqrt {\frac {4} {5}} E_n \leq \| t \|} C (\alpha, k) (1+ \| t \| ^ {| \alpha |}) \exp \Big (- \frac {3} {8} \| t \| ^ 2 \Big ) dt
$$
$$
\leq \int \Big (\| t \| ^ {- 1} \Big (\frac {4} {5} \Big) ^ {\frac {1} {2}} c_1 (k, k + 3) \min \Big\{\eta_ {k + 3} ^ {- \frac {1} {k + 3}}, \eta_ {k + 3} ^ {- \frac {1} {k + 1}}\Big \} \Big) ^ {- k-3} 
$$
$$\times C (\alpha, k) (1+ \| t \| ^ {| \alpha |}) \exp \Big (- \frac {3} {8} \| t \| ^ 2 \Big) dt
\leq C_4 (\alpha, k) \delta ^ 4_ \theta $$
and
$$
I_5 = \int \Big | D ^ \alpha \sum_ {r = 1} ^ {k} P_r (it, \kappa_ \nu) \exp \Big (- \frac {1} {2} \langle Dt, t \rangle \Big) \Big | dt
$$
$$
\leq \sum_ {r = 2} ^ {k} C (\alpha, k, r) \rho_ {r + 2}(1+\rho_{2}^{r-1}) + C_6 (\alpha, k) \sum_ {| \nu | = 3} | \kappa_ \nu |
$$
$$
\leq C_5 (\alpha, k) \delta_ \theta ^ 4 + C_6 (\alpha, k) \sum_ {| \nu | = 3} \Big | \sum_ {j = 1} ^ n \theta_j ^ 3 \kappa_ {\nu} (Z_j) \Big |
$$
$$
= C_5 (\alpha, k) \delta_ \theta ^ 4 + C_6 (\alpha, k) \sum_ {| \nu | = 3} \Big | \sum_ {j = 1} ^ n \theta_j ^ 3 \mu_ {\nu} (Z_j) \Big |.
$$
We get that there is a subset of weight coefficients $ \mathbb {Q} _1 \cap \mathbb {Q} _2
$
with  measure
$$
\lambda_{n-1} (\mathbb {Q} _1 \cap \mathbb {Q} _2) \geq1- \grave {C} _ {2} (k) \exp \Big (- \grave {c} _ {2} (k) \frac {n} {\delta ^ 4} \Big) - \grave {C} _1 \exp (- \grave {c} _1n) $$
$$ \geq 1- \grave {C} _ {5} (k) \exp \Big (- \grave {c} _ {5} (k) \frac {n} {\delta ^ 4} \Big) ,$$
such that for any vector of weight coefficients $ (\theta_1, \theta_2 ,\dots, \theta_n) \in \mathbb {Q} _1 \cap \mathbb {Q} _2 $
one has
\begin{equation}\label{temp_eq}
\Big | \int \limits_ {B} d (F_X- \Phi) \Big | \leq C_7 (k)\Big( \delta_ \theta ^ 4 + \frac {\delta ^ 4} {n} +  \sum_ {| \nu | = 3} \Big | \sum_ {j = 1} ^ n \theta_j ^ 3 \mu_ {\nu} (Z_j) \Big |\Big).
\end{equation}
Note that if $ \delta_ \theta ^ 4> (8k)^{-1}, $  the inequality (\ref{temp_eq}) holds automatically for a certain choice of universal constants.  Also, without loss of generality, we can require that
$$ \log ^ 2 \Big (\frac {\rho} {2 \grave {C} _ {5} (k)} \Big) \leq \grave {c} _ {5} (k) \frac { n} {\delta ^ 4}, $$
similarly, otherwise the statement of the Theorem holds for a special choice of the constant in the inequality. Further,
$$ \frac {\rho} {2 \grave {C} _ {5} (k)} \geq \exp \Big (- \grave {c} _ {5} (k) \Big (\frac {n } {\delta ^ 4} \Big) ^ {\frac {1} {2}} \Big)
$$
and
$$
\lambda_{n-1}( \mathbb {Q}_1 \cap \mathbb {Q}_2)> 1- \grave {C}_{5} (k) \frac {\rho} {2 \grave {C}_ {5} (k)} \geq 1- \frac {\rho} {2}.
$$
By Lemma \ref {lemma7} there exists a subset $ \mathbb {Q} _3 $ with  $\lambda_{n-1}(\mathbb {Q} _3) \geq 1-\rho/2, $
for which
$$  \delta_ \theta ^ 4 +  \sum_ {| \nu | = 3} \Big | \sum_ {j = 1} ^ n \theta_j ^ 3 \mu_ {\nu} (Z_j) \Big |  \leq \hat {C} _7 \Big (\log \Big (\frac {1} {\rho} \Big) \Big) ^ 2 \frac {\delta ^ 4} {n} +
\hat {C} _8 (k) \Big (\log \Big (\frac {1} {\rho} \Big) \Big) ^ {\frac {3} {2}} \frac {\delta ^ 4 } {n}.$$
Finally, for any vector of coefficients $ (\theta_1, \theta_2 ,\dots, \theta_n) \in \; \mathbb {Q} = \mathbb {Q} _1 \cap \mathbb {Q} _2 \cap \mathbb { Q} _3 $ we have
$$
   \sup \limits_ {B \in \mathfrak {B}} \Big | \int \limits_ {B} d (F_X- \Phi) \Big | \leq \Big (\log \Big (\frac {1} {\rho} \Big) \Big) ^ 2C (k) \frac {\delta ^ 4} {n} \leq C (\rho, k) \frac {\delta ^ 4} {n},
$$
moreover $$ \lambda_{n-1} (\mathbb {Q}) \geq 1- \frac {\rho} {2} - \frac {\rho} {2} = 1- \rho. $$
 
\end {proof}


\begin {thebibliography}{3}
\bibitem{S68} V. V. Sazonov, On the multidimensional central limit theorem, Sankhya Ser. A. 30 (1968), 181--204. 

\bibitem {KS}
B. Klartag, S. Sodin, Variations on the Berry--Esseen theorem. - Theory of Probability and its Applications, 2011, v. 56, no. 3, p. 514-533.

\bibitem{B20} S. G. Bobkov, Edgeworth Corrections in Randomized Central Limit Theorems. In: Klartag B., Milman E. (eds) Geometric Aspects of Functional Analysis. Lecture Notes in Mathematics, 2256 (2020). Springer, Cham. 71--97.

\bibitem{BCG18} S. G. Bobkov, G. P. 
Chistyakov, F. G\"otze, Berry--Esseen bounds for typical weighted sums. Electron. J. Probab. 23 (2018), no. 92, 1--22.

\bibitem{GNU17} F. G\"otze, A. A. Naumov, V. V. Ulyanov, Asymptotic Analysis of Symmetric Functions. J Theor Probab 30 (2017), 876--897. 

\bibitem{12} C.-G.~Esseen, Fourier analysis of distribution functions. A~mathematical study of the Laplace--Gaussian law. Acta Math. 77 (1945),  1--125.
 
\bibitem{20} Yu.~V.~Prokhorov, V.~V.~Ulyanov,
 Some approximation problems in statistics and probability. In 
<< Limit theorems in probability, statistics and number theory>>,  Springer Proc. Math. Stat. vol. 42,  Springer Heidelberg (2013), 235--249.

\bibitem{14} F.~G\"{o}tze, A.~Yu.~Zaitsev,
 Explicit rates of approximation in the CLT for quadratic forms.  Ann. Probab. 42(1) (2014), 
 354--397.

\bibitem{23} V.~V.~Ulyanov, Asymptotic Expansions for Distributions of Sums of Independent Random Variables in~$H$. Theory Probab. Appl. 31(1) (1987), 25--39.

\bibitem {B}
R. Bhattacharya, R Rao Ranga, Normal Approximation and Asymptotic Expansions. - New York, Wiley, 1976.
\bibitem {Saz}
V. Sazonov, Normal Approximation - Some Recent Advances. - New York, Springer-Verlag Berlin Heidelberg, 1981.
\bibitem {A}
R. Adamczak, A. Litvak, A. Pajor, Tomczak-Jaegermann N. Restricted isometry
property of matrices with independent columns and neighborly polytopes by random
sampling. - Constr. Approx., 2011, v. 34, no. 1, \\p. 61--88.
\bibitem {U}
6. S. Barsov, V. Ulyanov, Difference of Gaussian measures. - Journal of Soviet Mathematics, 1987, v. 5, no. 38, p. 2191-2198.
\bibitem {S}
E. Carlen, E. Lieb, M. Loss, A sharp analog of Young's inequality on SN and relatedentropy inequalities. - J. Geom. Anal, 2004, v. 4, no. 3, p. 487-520.

%\bibitem{KS11}  B. Klartag, S. Sodin, Variations on the Berry-Esseen theorem. (Russian) Teor. Veroyatn. Primen. 56 (2011), no. 3, 514--533; translation in Theory Probab. Appl. 56 (2012), no. 3, 403--419.
%
%\bibitem{BR} R.N. Bhattacharya, R. Ranga Rao, Normal approximation and asymptotic expansions. SIAM: Philadelphia, 2010.
% 
%\bibitem{S81} V. V. Sazonov, Normal Approximation -- Some Recent Advances. Springer--Verlag: Berlin, vol.879, 1981.


\end {thebibliography}
\end {document}